\documentclass[leqno]{amsart}
\input{diagxy}
\let\BarrSquare\square
\usepackage{amssymb}
\usepackage{stmaryrd}
\newtheorem{theorem}{Theorem}[section]
\newtheorem{lemma}[theorem]{Lemma}
\newtheorem{corollary}[theorem]{Corollary}
\newtheorem{proposition}[theorem]{Proposition}
\theoremstyle{definition}
\newtheorem{definition}[theorem]{Definition}
\newtheorem{notation}[theorem]{Notation}
\newtheorem{example}[theorem]{Example}
\newtheorem{remark}[theorem]{Remark}
\theoremstyle{remark}
\numberwithin{equation}{section}

\newcommand{\cSet}{\mathbf{Set}}
\newcommand{\cCa}{\mathbf{Cal}}
\newcommand{\cCas}{\mathbf{Cal^0_{Sym}}}
\newcommand{\cCasm}{\mathbf{Cal_{Sym}}}
\newcommand{\clo}{\mathsf{Cl}}
\newcommand{\cC}{\mathbf{C}}
\newcommand{\A}{\mathsf{A}}
\newcommand{\B}{\mathsf{B}}
\newcommand{\C}{\mathsf{C}}
\newcommand{\D}{\mathsf{D}}
\renewcommand{\l}{\mathcal{L}}
\newcommand{\cChu}{\mathbf{Chu}}
\renewcommand{\smash}{\mathbin{\sharp}}
\newcommand{\op}{{\mathrm{op}}}
\newcommand{\F}{\mathcal{F}}
\newcommand{\G}{\mathcal{G}}

\newcommand{\fk}{\mathcal{K}}
\newcommand{\lp}{\left(}
\newcommand{\rp}{\right)}
\newcommand{\lac}{\left\{}
\newcommand{\rac}{\right\}}
\newcommand{\lcr}{\left [}
\newcommand{\rcr}{\right ]}
\newcommand{\s}{\Sigma}
\newcommand{\aut}{\mathsf{Aut}}
\newcommand{\Sep}{\mathcal{S}}
\newcommand{\pro}{\mathsf{P}}
\newcommand{\p}{\perp}
\newcommand{\ot}{\otimes}
\newcommand{\otchu}{\otimes_{_C}}
\begin{document}
\title[The Chu construction for atomistic lattices]{The Chu construction for complete atomistic coatomistic lattices}
\author{Boris Ischi}
\address{Boris Ischi, Coll\`ege de Candolle, 5 rue
d'Italie, 1204 Geneva, Switzerland}
\email{boris.ischi@edu.ge.ch}
\thanks{The first author was partially supported by the Swiss National
Science Foundation.}
\author{Gavin J. Seal}
\address{Gavin J. Seal, \'Ecole Polytechnique F\'ed\'erale de Lausanne,
Switzerland}
\email{gavin\underbar{ }seal@fastmail.fm}
\begin{abstract}
The Chu construction is used to define a
$\ast$-autonomous structure on a category of complete atomistic coatomistic
lattices. This construction leads to a new tensor product
that is compared with a certain number of other existing
tensor products.
\end{abstract}
\maketitle
\section{Introduction}

In \cite{Shmuely:1974} Shmuely defines the tensor product of two partially ordered
sets as the set of all Galois connections between them
(see Definition \ref{DefinitionGaloisConnections}).
When the ordered sets are
complete lattices, this tensor product can equivalently be obtained by way of a
universal property in the category of complete lattices and join-preserving maps.
This tensor product appears as an instance of a result by Banaschewski and Nelson \cite{Banaschewski/Nelson:1976}
that is based on existence of an internal hom. Alternatively,
if no obvious internal hom is available, the Chu construction \cite{Barr:handbook}
proves to be useful to exhibit a suitable bifunctor.

In this work, we follow the second lead to obtain a tensor product $\circledast$ in a category of complete
{\it atomistic} lattices. Our original motivation for the study of this a
category is the Piron-Aerts approach to quantum logic, in which a physical system
is described by a complete atomistic lattice and time-evolution is modeled by a join-preserving
map sending atoms to atoms or zero \cite{Faure/Moore/Piron:1995}.

In order to use the Chu construction, we first need to restrict attention
to a suitable subcategory $\cCas$ of lattices, as discussed in Section \ref{setting}.
The rest of this paper is organized as follows.
In Section
\ref{SectionTheCategoryChu(Set_0,2_0)}, we recall the
definition of
$\ast$-autonomous categories and Chu spaces.
The category $\cCas$ and the bifunctor $\circledast$ are
introduced in Sections \ref{SectionTheCategoryCalSym^0} and
\ref{SectionTheBifunctor*} respectively. Section
\ref{Section*-autonomousStructureOnCalSym^0} is devoted to our
main result, namely $\cCas$ with $-\circledast-$ is
$\ast$-au\-to\-no\-mous.
In Section \ref{SectionComparison}, the tensor
product $-\circledast-$ is compared to other lattice-theoretical
tensor products. It is characterized in terms of a universal
property with respect to what we call {\it weak bimorphisms} in
Section \ref{SectionWeakBimorphisms}. Finally, in Section \ref{SectionDAC-lattices}
we focus on DAC-lattices as these play a central role in quantum logic.
For lattice-theoretic terminology, we refer to Maeda and Maeda
\cite{Maeda/Maeda:handbook} and for general terminology
concerning category theory, we refer to Mac Lane
\cite{McLane:handbook}.
\section{Setting}\label{setting}

As mentioned in the introduction, the category of complete atomistic
lattices with join-preserving maps sending atoms to atoms or $0$ appears naturally
in the foundations of quantum theory \cite{Faure/Moore/Piron:1995}
(see~\cite{Coecke/Moore/Wilce:handbook} for further references). The problem of describing compound systems in quantum axiomatic is related to the question of finding
a suitable tensor product for their models (see for instance \cite{Aerts/Daubechies:1979}, \cite{Valckenborgh:2008} or
\cite{Aerts02thelinearity} and references therein). Our goal is to apply the Chu construction to a subcategory of the category mentioned above
in order to find a new tensor product. The starting point is the trivial remark that a complete atomistic lattice, seen as a simple closure space, is a Chu space and that a join-preserving map sending  atoms to atoms defines an arrow in the Chu category. In this section, we outline the technical reasons that brought us
to study the subcategory $\cCas$ defined in section \ref{SectionTheCategoryCalSym^0}.

Denote by $2$ the two-element set $\{0,1\}$. The objects of the category
$\cChu_2\equiv\cChu(\cSet,2)$ are triples $(A,r,X)$, where $A$
and $X$ are sets, and $r$ is a map; $r:A\times X\rightarrow 2$. Arrows
are pair of maps $(f,g):(A,r,X)\rightarrow (B,s,Y)$ with
$f:A\rightarrow B$ and $g:Y\rightarrow X$ such that
$s(f(a),y)=r(a,g(y))$ for all $a\in A$ and $y\in Y$. The functor
$-^\bot$ is defined on objects as $(A,r,X)^\bot=(X,\check{r},A)$,
with $\check{r}(x,a)=r(a,x)$, and the bifunctor $-\ot-$ as
$\A_1\ot\A_2=(A_1\times A_2,t,\cChu_2(\A_1,\A_2^\bot))$ where
$\A_i=(A_i,r_i,X_i)$ and $t((a_1,a_2),(f,g))=r_2(a_2,f(a_1))$.

Obviously, if $\l$ is an atomistic lattice, then
$\G(\l):=(\s,r,\l)$, where $\s$ denotes the set of atoms of $\l$
and $r(p,a)=1\iff p\leq a$, is an object of $\cChu_2$.
If in addition of being atomistic lattices,
$\l_1$ and $\l_2$ are moreover complete, then
$(f,g)\in\cChu_2(\G(\l_1),\G(\l_2))$ if and only if there is a
Galois connection
(see Definition \ref{DefinitionGaloisConnections})
$(h,h^\circ)$ between $\l_1$ and $\l_2$ with $h$
sending atoms to atoms, $f=h\vert_{_{\s_1}}$ and $g=h^\circ$. As a
consequence, for the category $\cCa$ of complete atomistic
lattices with maps preserving arbitrary joins and sending atoms to
atoms, we have a full and faithful functor
$\G:\cCa\rightarrow\cChu_2$. Moreover, it can be proved that
$\cCa$ is closed under the tensor product of $\cChu_2$, {\it i.e.}
$\G(\l_1\varovee\l_2)\cong\G(\l_1)\ot\G(\l_2)$ where $\varovee$
denotes the complete lattice tensor product of Shmuely
\cite{Shmuely:1974} (or the tensor product of Fraser
\cite{Fraser:1976}, see Proposition \ref{PropositionFunctorG},
Remark \ref{RemarkShmuely} and \ref{RemarkFraser} below). However,
$\cCa$ is obviously not closed under $-^\bot$.

A natural way to get closure under $-^\bot$, is to replace in the
definition of the functor $\G$ the set $\l$ by the set of coatoms
of $\l$. Indeed, for the (non-full) subcategory $\cCasm$ of
coatomistic lattices with arrows $f$ having a right adjoint
$f^\circ$ sending coatoms to coatoms, we can define a full and
faithful functor $\fk:\cCasm\rightarrow\cChu_2$ as
$\fk(\l)=(\s,r,\s')$ and $\fk(f)=(f,f^\circ)$, where $\s'$ stands
for the set of coatoms of $\l$. Then $\fk(\l^\op)=\fk(\l)^\bot$
(where $\l^\op$ denotes the dual of $\l$ defined by the converse
order relation) but the category $\cCasm$ is not closed under the
tensor product of $\cChu_2$. For instance, consider $\mathsf{MO}_n$ the
complete atomistic lattice with $n$ atoms such that $1$ covers
each atom. Let $\l_1=\mathsf{MO}_n$ and $\l_2=\mathsf{MO}_m$ with
$n\ne m$. Then it is easy to check that
$\cChu_2(\fk(\l_1),\fk(\l_2)^\bot)=\emptyset$ since there is no
bijection between $\s_1$ and $\s_2$. As a consequence, there is
no $\l\in\cCasm$ such that
$\fk(\l)\cong\fk(\l_1)\ot\fk(\l_2)$.

In the preceding example, the reason why $\cCasm$ is not closed
under the tensor product of $\cChu_2$ is that
$\cCasm(\l_1,\l_2^\op)=\emptyset$. A simple way to remedy to this
is to consider more morphisms, namely maps preserving arbitrary
joins and sending atoms to atoms or $0$. Call $\cCa^0$ the
category of complete atomistic lattices equipped with those
morphisms. Then $\cCa^0$ can be embedded canonically in the
category $\cChu_{2_0}\equiv \cChu(\cSet_0,2_0)$, where $\cSet_0$
denotes the category of pointed sets with monoidal structure given
by the smash product, and $2_0$ denotes the set $\{0,1\}$ pointed
by $0$ (see Section \ref{SectionTheCategoryChu(Set_0,2_0)} for details). Indeed, it is easy to check that the functor
$\G^0:\cCa^0\rightarrow\cChu_{2_0}$ defined on objects as
$\G^0(\l)=(\s\cup\{0\},r,\l)$ where $\s\cup\{0\}$ is pointed
by $0$ and $\l$ is pointed by $1$, and where
$r(x,a)=0\iff x\leq a$, and on morphisms as
$\G^0(f)=(f,f^\circ)$, is full and faithful, and that $\cCa^0$ is
closed under the tensor product of $\cChu_{2_0}$.

However, as $\cCa$, the category $\cCa^0$ is not closed under the
dualizing functor $-^\bot$. Again, a natural way to obtain
closure under $-^\bot$ is to replace in the definition of $\G^0$
the set $\l$ by $\s'\cup\{1\}$
(we refer the reader to Section \ref{SectionTheCategoryCalSym^0} for the definition of $\cCas$).
Hence, we define a functor
$\F:\cCas\rightarrow\cChu_{2_0}$ as
$\F(\l)=(\s\cup\{0\},r,\s'\cup\{1\})$ and
$\F(f)=(f,f^\circ)$.
Then, for $\F$ to be full we have to
consider more morphisms than in $\cCasm$, namely maps preserving
arbitrary joins and sending atoms to atoms or $0$ with right
adjoint sending coatoms to coatoms or $1$.

In order to check closure under tensor product, let $\l_1$ and $\l_2$ be
complete atomistic coatomistic lattices. Then
\[
\F(\l_1)\otchu\F(\l_2)\cong
(\s_1\times\s_2\cup\{0\},t,\cChu_{2_0}(\F(\l_1),\F(\l_2)^\bot)),
\]
where $\otchu$ denotes the tensor product of $\cChu_{2_0}$. Now, by definition
$(f,g)\in\cChu_{2_0}(\F(\l_1),\F(\l_2)^\bot)$ if and only if
$f:\s_1\rightarrow \s_2'\cup\{1\}$, $g:\s_2\rightarrow
\s_1'\cup\{1\}$, and $q\leq f(p)\iff p\leq g(q)$,
for all atoms $p\in\s_1$ and $q\in\s_2$. To the map $f$ we can
associate a subset $x^f$ of $\s_1\times\s_2$ defined as
$x^f=\bigcup\{\{p\}\times\s_2[f(p)]\, ;\, p\in\s_1\}$, where $\s_2[b]$
denotes the set of atoms under $b$. Hence, it can be seen that
$\F(\l_1)\otchu\F(\l_2)\cong (\s_1\times\s_2\cup\{0\},
r,
\Gamma)$,
where $\Gamma$ is the set of all subsets $x^f$ of $\s_1\times\s_2$
and where for $X\in\Gamma$ and $p\in\s_1\times\s_2$, $r(0,X)=1$ and
$r(p,X)=1\iff p\in X$.
As
a consequence, $\F(\l_1)\otchu\F(\l_2)$ is in the image of $\F$
only if there is $\l\in\cCas$ such that $\s=\s_1\times\s_2$ and
$\{\s[x]\,;\,x\in\s'\cup\{1\}\}=\Gamma$, that is only if for
any $a,\, b\in\Gamma$ different from $\s_1\times\s_2$, $a$ is not a
subset of $b$ and $b$ is not a subset of $a$. This fails to be
true, for instance if $\l_1$ and $\l_2$ are powerset lattices (see
Example \ref{ExamplePowerSets}).

Therefore, in order to have closure under tensor product, the objects in
$\cCas$ cannot be all complete atomistic coatomistic lattices, but
me must impose some condition. We will prove that a sufficient
condition (which we call $\mathbf{A_0}$) is to ask that for any
two atoms $p$ and $q$ and any two coatoms $x$ and $y$, there is a
coatom $z$ and an atom $r$ such that $p\wedge z=0=q\wedge z$
and $r\wedge x=0=r\wedge y$. Note that our Axiom
$\mathbf{A_0}$ implies that the lattices are irreducible. We will
give an example of a complete atomistic orthocomplemented lattice
$\l$ which is irreducible but does not satisfy $\mathbf{A_0}$, and
such that there is no $\l_0\in\cCas$ with
$\F(\l_0)\cong\F(\l)\otchu\F(\l)$ (see Example
\ref{ExampleA0Necessary}).

Using the functor $\F$ we prove that $\cCas$ is closed under both
$-^\bot$ and the tensor product of $\cChu_{2_0}$, hence that $\cCas$
inherits the
$\ast$-autonomous structure of $\cChu_{2_0}$. This
result is presented in Theorem \ref{TheoremTheTheorem}, Section
\ref{Section*-autonomousStructureOnCalSym^0}.
\section{The category
$\cChu(\cSet_0,2_0)$}\label{SectionTheCategoryChu(Set_0,2_0)}

We begin by briefly recalling the definition of a
$\ast$-autonomous category. For details, we refer to Barr
\cite{Barr:handbook}, \cite{Barr:1991}.

\begin{definition} An autonomous category $\mathbf{C}$ is a monoidal
symmetric closed category. {\it Monoidal symmetric} means that
there is a bifunctor $-\ot-:{\bf C}\times{\bf C}\rightarrow{\bf C}$,
an object $\top$, and natural isomorphisms
$\alpha_{_{\A\B\C}}:(\A\ot \B)\ot \C\rightarrow \A\ot(\B\ot \C)$,
$r_{_{\A}}:\A\ot\top\rightarrow \A$, $l_{_{\A}}:\top\ot
\A\rightarrow \A$, and $s_{_{\A\B}}:\A\ot \B\rightarrow \B\ot \A$,
satisfying some coherence conditions (see the appendix). {\it
Closed} means that there is a bifunctor $-\multimap-:{\bf
C}^{\op}\times{\bf C}\rightarrow{\bf C}$ such that for all objects
$\A,\,\B,\,\C$ of ${\bf C}$, there is an isomorphism ${\bf
C}(\A\ot \B,\C)\cong{\bf C}(\B,\A\multimap \C)$, natural in $\B$
and $\C$.
\end{definition}

\begin{remark}
Let ${\bf C}$ be an autonomous category and $\bot$ an object of
${\bf C}$. Since ${\bf C}$ is closed and symmetric, for each
object $A$ we have
\[\begin{split}{\bf C}(\A\multimap\bot,\A\multimap\bot)&\cong {\bf
C}(\A\ot(\A\multimap\bot),\bot)\\
&\cong{\bf C}((\A\multimap\bot)\ot \A,\bot)\cong{\bf
C}(\A,(\A\multimap\bot)\multimap\bot).\end{split}\]
Hence, to the identity arrow $(\A\multimap\bot)\rightarrow
(\A\multimap\bot)$ corresponds an arrow
$\A\rightarrow((\A\multimap\bot)\multimap\bot)$.
\end{remark}

\begin{definition}[Barr, \cite{Barr:handbook}] If for every object
$\A$ of ${\bf C}$ the aforementioned arrow
$\A\rightarrow((\A\multimap \bot)\multimap\bot)$ is an
isomorphism, the object $\bot$ is called a {\it dualizing object}.
A $\ast$-autonomous category ${\bf C}$ is an autonomous category
with a dualizing object. Usually, $\A\multimap\bot$ is written
$\A^\bot$.
\end{definition}

Chu's paper \cite{Barr:handbook} (see also Barr  \cite{Barr:1996},
\cite{Barr:1999}) describes a construction of a
$\ast$-auto\-no\-mous category starting with a finitely complete
autonomous category. We outline the construction of Chu for the
category $\cSet_0$ of pointed sets and pointed maps.

\begin{definition}
On the category $\cSet_0$ of pointed sets with pointed maps, we
define the {\it smash product}
$-\smash-:\cSet_0\times\cSet_0\rightarrow\cSet_0$ as
\[A\smash B:=[(A\backslash\{0_A\})\times (B\backslash\{0_B\})]\cup
\{0_\sharp\},\]
where $0_A$ and $0_B$ are the respective base-points of $A$ and
$B$. Moreover, we write $2_0$ for the set $\{0,1\}$ pointed by
$0$.
\end{definition}

\begin{lemma}\label{LemmaPointedSetsisMonoidalClosed}
The category $\langle\cSet_0,\smash,2_0,\multimap\rangle$, with
$A\multimap  B:=\cSet_0(A, B)$ pointed by the constant map, is
autonomous.
\end{lemma}

\begin{proof} The proof is direct and is omitted here.
\end{proof}

\begin{notation} Let $A,\, B$ and $C$ be pointed sets and let
$r:A\smash B\rightarrow C$ be a pointed map. We do not distinguish
this map from the map defined on $A\times B$ with value in $C$
such that $r(a,b)=0_C$, whenever $a=0_A$ or $b=0_B$.
\end{notation}

\begin{definition}
An object of $\cChu_{2_0}:=\cChu(\cSet_0,2_0)$ is a triplet
$(A,r,X)$, where $A$ and $X$ are pointed sets, and $r:A\smash
X\rightarrow 2_0$ is a pointed map. An arrow is a pair of pointed
maps $(f,g):(A,r,X)\rightarrow(B,s,Y)$, with $f:A\rightarrow B$
and $g:Y\rightarrow X$, satisfying $s(f(a),y)=r(a,g(y))$ for all
$a\in A$ and $y\in Y$.
\end{definition}

\begin{remark} Let $\A=(A,r,X)$ and $\B=(B,s,Y)$ be objects of
$\cChu_{2_0}$. The pair of constant maps $f:A\rightarrow
B;\,a\mapsto 0_B$ and $g:Y\rightarrow X;\,y\mapsto 0_X$ forms an
arrow of $\cChu_{2_0}(\A,\B)$ which we call the {\it constant}
arrow.\end{remark}

\begin{definition}\label{DefinitionTensorChu} For $i=1,2$, let
$\A_i=(A_i,r_i,X_i)$ and $\B_i=(B_i,s_i,Y_i)$ be objects of
$\cChu_{2_0}$, and $(f_i,g_i)\in \cChu_{2_0}(\A_i,\B_i)$.

The functor $-^\bot:\cChu_{2_0}^\op\rightarrow\cChu_{2_0}$ is
defined on objects as $\A_1^\bot:=(X_1, \check{r}_1,A_1)$, where
$\check{r}_1(x,a)=r_1(a,x)$, and on arrows as
$(f_1,g_1)^\bot:=(g_1,f_1):\B_1^\bot\rightarrow\A_1^\bot$.

The bifunctor
$-\otchu-:\cChu_{2_0}\times\cChu_{2_0}\rightarrow\cChu_{2_0}$ is
defined on objects as
\[\A_1\otchu\A_2=(A_1\smash A_2,t,\cChu_{2_0}(\A_1,\A_2^\bot)),
\]
with $\cChu_{2_0}(\A_1,\A_2^\bot)$ pointed by the constant arrow,
and with $t$ defined as
\[t((a_1,a_2),(f,g)):=r_1(a_1,g(a_2))=\check{r}_2(f(a_1),a_2).\]
Further, $f_1\otchu f_2:\A_1\otchu\A_2\rightarrow\B_1\otchu\B_2$
is defined as $(a_1,a_2)\mapsto (f_1(a_1),f_2(a_2))$ if
$f_1(a_1)\ne 0$ and $f_2(a_2)\ne 0$ (with $(f_1\otchu
f_2)(a_1,a_2):=0_\sharp$ if $f_1(a_1)=0$ or $f_2(a_2)=0$), and as
$(f,g)\mapsto(g_2\circ f\circ f_1, g_1\circ g\circ f_2)$ for
$(f,g)\in\cChu_{2_0}(\B_1,\B_2^\bot)$.
\end{definition}

\begin{definition}\label{DefinitionHomSetChu} The object $\top$ is
defined as $\top:=(2_0, r,2_0)$ with $r$
the isomorphism between $2_0\smash 2_0$ and $2_0$.
Moreover, the
dualizing object $\bot$ is defined as $\top^\bot$. Finally, the
bifunctor $-\multimap-$ is given by $\A\multimap \B:=(\A\otchu
\B^\bot)^\bot$.\end{definition}

\begin{remark} The object $\top$ is the tensor product unit; $\A\otchu\top\cong
\A$. Hence, we have $\A\multimap \top^\bot\cong \A^\bot$, and
$\bot$ is the dualizing object.
Note that $\top\cong\bot$.
\end{remark}

Since $\cSet_0$ is finitely complete and autonomous (Lemma
\ref{LemmaPointedSetsisMonoidalClosed}), we have the following
result.

\begin{proposition} \label{LemmaChuPointedSetsis*autonomous}
The category $\langle \cChu_{2_0},
\otchu,\top,\multimap,\bot\rangle$ is $\ast$-autonomous.
\end{proposition}
\section{The category $\cCas$}\label{SectionTheCategoryCalSym^0}

\begin{definition}
Let $\s$ be a nonempty set and $\l\subseteq
2^\s$. We say that $\l$ is a {\it simple closure space} on $\s$ if
$\l$ contains $\emptyset$ and all singletons, and if $\l$
is closed under arbitrary set-intersections ({\it i.e.}
$\bigcap\omega\in\l$ for all $\omega\subseteq\l$),
in particular $\l$ contains $\s=\bigcap\emptyset$.
Note that a
simple closure space (ordered by set-inclusion) is a complete
atomistic lattice.
Recall that an element $p$ of a lattice with 0 is called an {\it atom} if $p$
{\it covers} $0$ ({\it i.e.} $0\leq x\leq p$ $\Rightarrow$ $x=0$ or $x=p$). Moreover, a lattice with 0 is
called atomistic when every non-zero element $a\in\l$ is the join of the atoms under $a$.
For $p\in\s$, we identify $p$ with
$\{p\}\in\l$.
\end{definition}

\begin{notation}
Let $\l$ be a poset and $a\in\l$. The bottom and top elements
of $\l$, if they exist, are denoted by $0$ and $1$ respectively.
We denote by $\l^\op$ the dual of $\l$ (defined by the
converse order relation), by $\s$ and $\s'$ the sets of atoms
and coatoms of $\l$ respectively, by $\s[a]$ the set of atoms
under $a$, and by $\s'[a]$ the set of coatoms above $a$.
For any subset $\omega\subseteq\l$, we write
\[\clo(\omega)=\{\s[a]\,;\,a\in\omega\}\subseteq 2^{\s},\]
ordered by set-inclusion. Note that if $\l$
is a complete atomistic lattice, then $\clo(\l)$ is a simple
closure space on the set of atoms of $\l$.
\end{notation}

\begin{definition}\label{DefinitionGaloisConnections}
Let $\l_1$ and $\l_2$ be posets. A Galois connection between
$\l_1$ and $\l_2$ (or equivalently an adjunction) is a pair
$(f,g)$ of order-preserving maps with $f:\l_1\rightarrow\l_2$ and
$g:\l_2\rightarrow\l_1$ such that for any $a\in\l_1$ and
$b\in\l_2$, $f(a)\leq b\iff a\leq g(b)$.

Let $\l_1$ and $\l_2$ be complete lattices and
$f:\l_1\rightarrow\l_2$ a map. The join and meet in $\l_i^\op$ are
denoted by $\bigvee^\op$ and $\bigwedge^\op$ respectively. The map
$f^\circ:\l_2\rightarrow\l_1$ is defined as
\[f^\circ(b):=\bigvee\{a\in\l_1\,;\,f(a)\leq b\},\]
and $f^\op:\l_2^\op\rightarrow\l_1^\op$ as $f^\op(b):=f^\circ(b)$.
Finally, $2$ stands for the lattice with only two elements.
\end{definition}

\begin{remark}
Let $\l_1$ and $\l_2$ be posets and $(f,g)$ a Galois connection
between $\l_1$ and $\l_2$. Then $g=f^\circ$
preserves all existing meets and both $f$ and $f^\op$ preserve all existing joins.
Moreover, when $\l_1$ and $\l_2$ are complete lattices and $h:\l_1\rightarrow \l_2$
is a map then $h$ preserves arbitrary joins $\iff$
$(h,h^\circ)$ is a Galois connection between $\l_1$ and $\l_2$.
\end{remark}

\begin{lemma}\label{LemmaJoin-preserving}
Let $\l_1$ and $\l_2$ be complete atomistic lattices and let
$f:\l_1\rightarrow \l_2$ be a map that
sends atoms to atoms or $0$. Denote by $F$ the restriction of $f$
to atoms. Then
$f$ preserves arbitrary joins $\iff$
$f(a)=\bigvee F(\s[a])$ and
$F^{-1}(\s[b]\cup\{0\})\in\clo(\l_1)$, $\forall a\in\l_1$,
$b\in\l_2$.
\end{lemma}

\begin{proof}
Suppose that $f$ preserves arbitrary joins. Let $p$ be an
atom under $\bigvee F^{-1}(\s[b]$ $\cup\{0\})$. Then
\[
f(p)\leq \bigvee F(F^{-1}(\s[b]\cup\{0\}))\leq b,
\]
hence $p\in F^{-1}(\s[b]\cup\{0\})$.

We now prove the converse. Define $g(b)=\bigvee F^{-1}(\s[b]\cup \{0\})$.
We prove that the pair $(f,g)$ forms a Galois
connection between $\l_1$ and $\l_2$
(see Definition \ref{DefinitionGaloisConnections})
Let $a\in\l_1$ and
$b\in\l_2$. Suppose that $f(a)\leq b$. Then
$\bigvee F(\s[a])\leq b$, hence,
$f(p)\leq b$ for all
$p\in\s[a]$. Therefore $a\leq g(b)$. Suppose now that $a\leq
g(b)$. Then for any $p\in\s[a]$, $p\leq g(b)$, hence, from the
second hypothesis, $f(p)\leq b$. As a consequence, $f(a)=\bigvee
F(\s[a])\leq b$.
\end{proof}

\begin{definition}
We denote by $\cCas$ the following category: the objects are all
complete atomistic coatomistic lattices $\l$ such that
\renewcommand\theequation{$\mathbf{A_0}$}
\begin{equation}
\forall\, x\in\s',\ \bigvee\lp\s\backslash\s[x]\rp=1\ \mbox{and}\ \forall
p\in\s,\ \bigwedge\lp\s'\backslash\s'[p]\rp=0.
\end{equation}
The
arrows are all maps $f$ preserving arbitrary joins and sending
atoms to atoms or $0$ such that $f^\op$ sends atoms to atoms or
$0$ ({\it i.e.} $f^\circ$ sends coatoms to coatoms or $1$).
\end{definition}

\begin{remark}
The condition $\forall\, x\in\s',\ \bigvee\lp\s\backslash\s[x]\rp=1$ is equivalent to
\renewcommand\theequation{$\mathbf{\dagger}$}
\begin{equation}\label{conditiondagger}
\forall\, x,\,y\in\s',\ \s[x]\cup \s[y]\ne \s
\end{equation}
so that axiom $\mathbf{(A_0)}$ is equivalent to requiring this condition for $\l$ and
$\l^\op$ (that is
$\forall\ p,\,q\in\s$, $\s'[p]\cup\s'[q]\ne \s'$.)
\end{remark}

\begin{remark}\label{RemarkOPfunctor}
Note that $2\in\cCas$. Moreover, the map $-^\op:\cCas^\op
\rightarrow\cCas$ is a functor. Indeed, consider two arrows of
$\cCas$, say $g:\l_1\rightarrow\l_2$ and $f:\l_2\rightarrow\l_3$.
Let $c\in\l_3^\op$. Then we have
\[\begin{split}g^\op\circ f^\op(c)&=g^\circ(f^\circ(c))=
\bigvee\{a\in\l_1\,;\,g(a)\leq
f^\circ(c)\}\\
&=\bigvee\{a\in\l_1\,;\,f(g(a))\leq c\}={(f\circ
g)}^\circ(c)={(f\circ g)}^\op(c).\end{split}\]
Finally, note that Axiom ${\mathbf A_0}$ will only be needed for
Lemma \ref{Lemma*circcoatomistic}; note also that it implies that
$\l$ is irreducible (see \cite{Maeda/Maeda:handbook}, Theorem
4.13).
\end{remark}

\begin{definition}
A lattice with $0$ has the {\it covering property} if for any atom
$p$ and any $a\in\l$, $p\wedge a=0$ implies that $p\vee a$
covers $a$ ({\it i.e.} $p\vee a\geq x\geq a$ $\Rightarrow$ $x=p\vee a$ or $x=a$).

A lattice $\l$ with $0$ and $1$ such that $\l$ and its dual $\l^\op$ are
atomistic with the covering property, is called a {\it
DAC-lattice}.
\end{definition}

\begin{example}
Let $\l$ be an irreducible complete DAC-lattice. Then $\l$ is an
object of $\cCas$.

Indeed, let $x,\,y$ be coatoms. Suppose that
$\s[x]\cup\s[y]=\s$. Let $z$ be a coatom above $x\wedge y$.
Let $p$ be an atom under $z$ such that $p\wedge x\wedge
y=0$. Since by hypothesis $\s[x]\cup\s[y]=\s$, one has $p\leq x$
or $p\leq y$. Note that since $\l^\op$ has the covering property,
$z$ covers $x\wedge y$. Hence, if $p\leq x$, then $z=p\vee
(x\wedge y)=x$, and if $p\leq y$, then $z=p\vee (x\wedge
y)=y$. As a consequence, the set of coatoms above $x\wedge y$
is given by $\{x,y\}$, a contradiction. Indeed, recall that the
join of any two atoms of an irreducible complete DAC-lattice
contains a third atom (see \cite{Maeda/Maeda:handbook}, Theorems
28.8 and 27.6, and Lemma 11.6), hence, by duality, for any two
coatoms $x$ and $y$, there is a third coatom above their meet.

Finally, by duality, $\l^\op$ also satisfies condition ${\mathbf (\dagger)}$ so that
$\l$ satisfies axiom ${\mathbf A_0}$.
\end{example}

We end this section by recalling the relation between irreducible
complete DAC-lattices and lattices of closed subspaces of vector
spaces.

\begin{definition}[see \cite{Maeda/Maeda:handbook}, Definition 33.1]
Let $E$ be a left vector space (respectively $F$ a right vector
space) over a division ring $\mathbb{K}$
and $f:E\times F\rightarrow\mathbb{K}$ a non-degenerate
bilinear map.
Then $(E,F,f)$ is called a
{\it pair of dual spaces}. For $A\subseteq E$,
define
\[A^\p:=\{y\in F\,;\, f(x,y)=0,\, \forall x\in A\},\]
and for $B\subseteq F$ define $B^\p$ similarly. Define
\[\l_F(E):=\{A\subseteq E\,;\,A^{\p\p}=A\},\]
ordered by set-inclusion.\end{definition}

\begin{theorem}[see \cite{Maeda/Maeda:handbook}, Theorem 33.4]
Let $(E,F)$ be a pair of dual spaces. Then $\l_F(E)$ is an
irreducible complete DAC-lattice.
\end{theorem}

\begin{theorem}[see \cite{Maeda/Maeda:handbook},
Theorem 33.7]\label{TheoremMaeda} If $\l$ is an irreducible
complete DAC-lattice of length $\geq 4$, then there exists a pair
of dual spaces $(E,F)$ over a division ring $\mathbb{K}$ such that
$\l\cong\l_F(E)$.\end{theorem}

\begin{remark}\label{RemarkFiniteSubspaces}
If $x$ is a finite dimensional subspace of $E$, then $x\in\l_F(E)$
(see \cite{Maeda/Maeda:handbook}, Lemma 33.3.2).\end{remark}
\section{The bifunctor $\circledast$}\label{SectionTheBifunctor*}

The bifunctor $\circledast$ will provide $\cCas$ with a suitable
``tensor product'' in order to make it into a
$\ast$-autonomous
category.

\begin{notation} For $p\in\s_1\times\s_2$, we denote  the first component of $p$ by
$p_1$ and the second by $p_2$ that is $p=(p_1,p_2)$. For $R\subseteq\s_1\times\s_2$, we
adopt the following notations.
\[\begin{split}
R_1[p]&:=\{q_1\in\s_1\,;\, (q_1,p_2)\in R\},\\
R_2[p]&:=\{q_2\in\s_2\,;\, (p_1,q_2)\in R\}.\end{split}\]
\end{notation}

\begin{remark}
Note that $R_1[p]$ (respectively $R_2[p]$) depends only on $p_2$
(respectively only on $p_1$). For $(r,s)\in\s_1\times\s_2$, we define
$R_1[s]$ as $R_1[(p,s)]$ and $R_2[r]$ as $R_2[(r,q)]$ for any
$(p,q)\in\s_1\times\s_2$.
\end{remark}

\begin{definition}\label{Definition*TensorProduct}
Let $\l_1$ and $\l_2$ be complete atomistic coatomistic lattices.
Then we define
\[\begin{split}
\s_\circledast':=\lac\vbox{\vspace{0.4cm}}
R\subsetneqq\s_1\times\s_2\,;\,\right.&\,R_1[p]\in
\clo(\s_1'\cup\{1\})\ \mbox{and}\\
&\left.R_2[p]\in\clo(\s_2'\cup\{1\}),\, \forall\,
p\in\s_1\times\s_2\rac
\end{split}\]
and
\[\l_1\circledast\l_2:=\lac\bigcap\omega\,;\,
\omega\subseteq\s_\circledast'\cup\{\s_1\times\s_2\}\rac,\]
ordered by set-inclusion.
\end{definition}

\begin{remark}\label{Remark2*tensorUnit}
Note that $\l_1\circledast 2\cong\l_1$.\end{remark}

\begin{notation} Let $\l_1$ and $\l_2$ be complete atomistic
lattices, $a_1\in\l_1$, and $a_2\in\l_2$. Then we define
\[\begin{split}
a_1\circ a_2&:=\s[a_1]\times\s[a_2],\\
a_1\mathbin{\square} a_2&:= (a_1\circ 1)\cup (1\circ
a_2)=(\s[a_1]\times\s_2)\cup(\s_1\times\s[a_2]).
\end{split}\]
\end{notation}

\begin{lemma}\label{LemmaCoatomsAertsCoatoms*} Let
$\l_1,\,\l_2\in\cCas$ and $\s_\varowedge':=\{x_1\mathbin{\square}
x_2\,;\,(x_1,x_2)\in\s_1'\times\s_2'\}$. Then
$\s_\varowedge'\subseteq\s_\circledast'$ and $\s_\varowedge'$ is a
set of coatoms of $\l_1\circledast\l_2$.\end{lemma}

\begin{proof}
Let $(x_1,x_2)\in\s_1'\times\s_2'$, $X=x_1\mathbin{\square} x_2$, and
$p\in\s_1\times\s_2$. Then $X_1[p]=\s_1$ if $p_2\leq x_2$ and
$X_1[p]=\s[x_1]$ otherwise. Similarly, $X_2[p]=\s_2$ if $p_1\leq x_1$
and $X_2[p]=\s[x_2]$ otherwise. As a consequence,
$X\in\s_\circledast'$, thus $X\in\l_1\circledast\l_2$.

Let $R\in\l_1\circledast\l_2$ such that $X\subseteq R$ and $X\ne
R$. Let $p\in R\backslash X$. Then $\{p_1\}\times (\{p_2\}\cup
\s[x_2])\subseteq R$, hence, since $x_2$ is a coatom of $\l_2$, by
Definition \ref{Definition*TensorProduct}, $\{p_1\}\times\s_2\subseteq
R$. As a consequence, $(\s[x_1]\cup\{p_1\})\times\s_2\subseteq R$ (and
$\s_1\times(\s[x_2]\cup\{p_2\})\subseteq R$). Hence, for all $s\in\s_2$,
$(\s[x_1]\cup\{p_1\})\times\{s\}\subseteq
R_1[s]\in\clo(\s_1'\cup\{1\})$, thus, since $x_1$ is a coatom of
$\l_1$, $R_1[s]=\s_1$ for all $s\in\s_2$, that is, $R=\s_1\times\s_2$.
As a consequence, $X$ is a coatom of $\l_1\circledast\l_2$.
\end{proof}

\begin{lemma}\label{Lemma*circcoatomistic}
Let $\l_1,\,\l_2\in\cCas$. Then $\l_1\circledast\l_2$ is a simple
closure space on $\s_1\times\s_2$. Moreover,
$\l_1\circledast\l_2\in\cCas$ and the set of coatoms of
$\l_1\circledast\l_2$ is given by $\s_\circledast'$.
\end{lemma}

\begin{proof} We first prove that $\l_1\circledast\l_2$ is a
simple closure space on $\s_1\times\s_2$. By Remark
\ref{Remark2*tensorUnit}, we can assume that $\l_1\ne 2$ and that
$\l_2\ne 2$. By definition, $\l_1\circledast\l_2$ contains
$\s_1\times\s_2$ and all set-intersections. Let $p\in\s_1\times\s_2$, and
\[\begin{split}X_1&:=\bigcap\lac x_1\mathbin{\square} x_2\,;\,
(x_1,x_2)\in\s'[p_1]\times\s'_2\rac,\\
X_2&:=\bigcap\lac x_1\mathbin{\square} x_2\,;\,
(x_1,x_2)\in\s'_1\times\s'[p_2]\rac.\end{split}\]
Then $p\in X_1\cap X_2$. By Lemma
\ref{LemmaCoatomsAertsCoatoms*},
$X_1,\,X_2\in\l_1\circledast\l_2$. Moreover
\[X_1=\bigcup\lac\lp\bigcap \lcr f^{-1}(1)\rcr_1\rp\times\lp\bigcap
\lcr f^{-1}(2)\rcr_2\rp\,;\,f\in 2^{\clo(\s'[p_1])\times\clo(\s'_2)}\rac,\]
where $\lcr f^{-1}(1)\rcr_1$ denotes the set of all $A\in \clo(\s'[p_1])$ such that there exists  $B\in \clo(\s'_2)$ with
$f(A,B)=1$.
Now, if $f^{-1}(1)\ne \clo(\s'[p_1])$, then $\bigcap
f^{-1}(2)=\emptyset$. Therefore, $X_1=p_1\circ 1$, and for the same
reason, $X_2=1\circ p_2$. Hence $\{p\}=X_1\cap X_2$, thus
$\{p\}\in\l_1\circledast\l_2$. As a consequence,
$\l_1\circledast\l_2$ is a simple closure space on $\s_1\times\s_2$.

To prove that the set of coatoms of $\l_1\circledast\l_2$ is given
by $\s_\circledast'$, it suffices to check that if $x,\
y\in\s_\circledast'\cup\{1\circ 1\}$ and $x\subsetneqq y$, then
$y=1\circ 1$. Let $p\in\s_1\times\s_2$ such that $x_2[p]\subsetneqq
y_2[p]$. Then $y_2[p]=\s_2$, since by definition, $\bigvee x_2[p]$ is
either a coatom of $\l_2$ or $1$. Let $q_2\not\in x_2[p]$. First,
since $x\subseteq y$, we have $x_1[q_2]\subseteq y_1 [q_2]$. Now,
by hypothesis, $p_1\not\in x_1[q_2]$ whereas $p_1$ is in $y_1
[q_2]$. As a consequence, $y_1 [q_2]=\s_1$, for any $q_2\not\in
x_2[p]$, therefore $\s_1\times (\s_2\backslash x_2[p])\subseteq y$. By
Axiom ${\mathbf A_0}$, it follows that $\bigvee(\s_2\backslash
x_2[p])=1$. Thus, we find that $y=1\circ 1$.

We now check that Axiom ${\mathbf A_0}$ holds in
$\l_1\circledast\l_2$. Let $x,\, y\in\s_\circledast'$ and
\[\begin{split}
&A:=\{p_1\in\s_1\,;\,x_2[p_1]=\s_2\}\\
&B:=\{p_1\in\s_1\,;\,y_2[p_1]=\s_2\}.
\end{split}\]
Note that
\[\begin{split}
A=\bigcap\{ x_1[s]\,;\,s\in\s_2\}\hspace{0.5cm}
\mbox{and}\hspace{0.5cm} B=\bigcap\{ y_1[s]\,;\,s\in\s_2\},
\end{split}\]
hence $A,\, B\in\clo(\l_1)$. Indeed,
\[\begin{split} r\in A&\iff x_2[r]=\s_2
\iff r\times\s_2\subseteq x\\ &\iff r\in x_1[s],\
\forall s\in\s_2\iff
r\in\bigcap\{x_1[s]\,;\,s\in\s_2\}.
\end{split}\]

Suppose now that $x\cup y=\s_1\times\s_2$. Then
for any $p\in\s_1\times\s_2$, $x_2[p]\cup y_2[p]=\s_2$. As a consequence,
since condition
${\mathbf (\dagger)}$ holds in $\l_2$, we have
that $x_2[p]\ne \s_2\Rightarrow y_2[p]=\s_2$, and $y_2[p]\ne
\s_2\Rightarrow x_2[p]=\s_2$; whence $A\cup B=\s_1$, a
contradiction, since condition ${\mathbf (\dagger)}$ holds in $\l_1$.

It remains to check that condition ${\mathbf (\dagger)}$ holds in
${(\l_1\circledast\l_2)}^\op$. Let $p,\,q\in\s_1\times\s_2$. Since
condition ${\mathbf (\dagger)}$ holds in $\l_1^\op$ and in $\l_2^\op$, there
is $(x_1,x_2)\in\s_1'\times\s_2'$ such that for $i=1$ and $i=2$,
$p_i\wedge x_i=q_i\wedge x_i=0$. As a consequence, we have
$p\wedge (x_1\mathbin{\square} x_2)=q\wedge (x_1\mathbin{\square} x_2)=0$.
Moreover, by Lemma \ref{LemmaCoatomsAertsCoatoms*}, $x_1\mathbin{\square}
x_2$ is a coatom of $\l_1\circledast\l_2$.
\end{proof}

\begin{lemma}\label{LemmaBijectionCoatoms*TensorCalSym(L1,L2^op)}
Let $\l_1,\,\l_2\in\cCas$. There is a bijection
\[\xi:\s_\circledast'\cup\{1_\circledast\}\rightarrow\cCas(\l_1,\l_2^\op),\]
such that for all $x\in\s_\circledast'\cup\{1_\circledast\}$,
we have $x=\{p\circ \xi(x)(p)\,;\,p\in\s_1\}$ (where $1_\circledast$ stands for $\s_1\times\s_2$).\end{lemma}

\begin{proof} Let $x\in \s_\circledast'\cup\{1_\circledast\}$.
Define $F_x:\s_1\cup\{0\}\rightarrow\s_2'\cup\{1\}$ as
\[F_x(p_1):=\bigvee x_2[p_1],\]
and $F_x(0)=1$. Moreover, define $f_x:\l_1\rightarrow\l_2^\op$ as
$f_x(a)=\bigvee^\op F_x(\s[a])$. Obviously, $f_x$ sends atoms to
atoms or $0$. Let $b\in\l_2^\op$ and
$A:=F_x^{-1}(\s'[b]\cup\{1\})$ ({\it i.e.}
$A=\{r\in\s_1\,;\,F_x(r)\leq^\op b\}$). Note that
\[\begin{split}
r\in A&\iff\s[b]\subseteq
x_2[r]\iff\{r\}\times\s[b]\subseteq x\iff
(r,s)\in x\, ,\,\forall s\in\s[b]\\
&\iff r\in x_1[s]\, ,\,\forall s\in\s[b]\iff
r\in \bigcap\{x_1[s]\,;\, s\in\s[b]\},
\end{split}\]
hence $A=\bigcap\{x_1[s]\,;\, s\in\s[b]\}$. As a consequence,
$A\in\clo(\l_1)$, therefore, by Lemma \ref{LemmaJoin-preserving},
$f_x$ preserves arbitrary joins.

Let $q$ be an atom of $\l_2$. Then
\[f_x^\circ(q)=\bigvee\{a\in\l_1\,;\,f_x(a)\leq^\op
q\}=\bigvee\{p\in\s_1\,;\, q\leq f_x(p)\}=\bigvee x_1[q].\]
Therefore, $f_x\in\cCas(\l_1,\l_2^\op)$.

Let $f\in\cCas(\l_1,\l_2^\op)$. Define $x^f\subseteq \s_1\times\s_2$
as $x^f=\bigcup\{p\circ f(p)\,;\,p\in\s_1\}$. Let
$p\in\s_1\times\s_2$. Then
\[x^f_2[p]=\s[f(p_1)]\in\clo(\s_2'\cup\{1\}),\]
and by Lemma \ref{LemmaJoin-preserving}
\[\begin{split}
x_1^f[p]&=\{r\in\s_1\,;\,p_2\leq
f(r)\}=\{r\in\s_1\,;\,f(r)\leq^\op
p_2\}\\
&=\s[\bigvee\{r\in\s_1\,;\,f(r)\leq^\op
p_2\}]=\s[\bigvee\{a\in\l_1\,;\,f(a)\leq^\op
p_2\}]\\
&=\s[f^\circ(p_2)]\in\clo(\s_1'\cup\{1\}).
\end{split}\]

Obviously, we have $x^{f_x}=x$ and $f_{x^f}=f$.\end{proof}

\begin{lemma}\label{LemmaAertsOtfBiFuntors}
For $i,\,j\in\{1,2\}$, let $\l_i^j\in\cCas$ and
$f_i\in\cCas(\l_i^1,\l_i^2)$. Then there is a unique
$u\in\cCas(\l_1^1\circledast\l_2^1,\l_1^2\circledast\l_2^2)$, with
$u(p)=(f_1(p_1),f_2(p_2))$, for any $p\in\s_1^1\times\s_2^1$ such that
$f_1(p_1)\ne 0$ and $f_2(p_2)\ne 0$. We denote $u$ by
$f_1\circledast f_2$.
\end{lemma}

\begin{proof} Write $F_i$ for $f_i$ restricted to atoms and define
$F$ as $F(p_1,p_2)=0_\circledast$ if $F_1(p_1)=0$ or $F_2(p_2)=0$,
and $F(p_1,p_2)=(F_1(p_1),F_2(p_2))$ otherwise.
Define $u$ as
\[
u(a)=\bigvee\{F(p_1,p_2)\,;\,(p_1,p_2)\in a\}.
\]
By Lemma
\ref{LemmaJoin-preserving}, it suffices to show that for any
coatom $x$ of $\l_1^2\circledast\l_2^2$,
$F^{-1}(x\cup\{0_\circledast\})$ is a coatom of
$\l_1^1\circledast\l_2^1$ or $1$. Let $p\in\s_1\times\s_2$. Suppose
that $f_2(p_2)\ne 0$. Then
\[\begin{split}
(F^{-1}(x\cup\{0_\circledast\}))_1[p]&=\{q_1\in\s_1\,;\,
(q_1,p_2)\in F^{-1}(x\cup\{0_\circledast\})\}\\
&=\{q_1\in\s_1\,;\, F(q_1,p_2)\in
x\cup\{0_\circledast\}\}\\
&=\{q_1\in\s_1\,;\,F_1(q_1)\ne 0\ \mbox{and}\ (F_1(q_1),F_2(p_2))\in
x\}\cup F_1^{-1}(0)\\
&=F_1^{-1}\lp\{r\in\s_1\,;\,(r,F_2(p_2))\in x\}\cup\{0\}\rp\\
&=F_1^{-1}(x_1[F_2(p_2)]\cup\{0\}).
\end{split}\]
If $f_2(p_2)=0$, then
$(F^{-1}(x\cup\{0_\circledast\}))_1[p]=\s_1$. As a consequence,
we find that for $i=1$ and $i=2$,
$(F^{-1}(x\cup\{0_\circledast\}))_i[p]\in\clo(\s_i'\cup\{1\})$
for any $p\in\s_1^1\times\s_2^1$.
\end{proof}

\begin{proposition}\label{Corollary*TensorisaBiFunctor}
$-\circledast-:\cCas\times\cCas\rightarrow\cCas$ is a
bifunctor.\end{proposition}
\section{
The $\ast$-autonomous structure on
$\cCas$}\label{Section*-autonomousStructureOnCalSym^0}

We can now turn to our main result. The functor given in the
following definition explains where the tensor product of $\cCas$
comes from. It is also useful to understand the
$\ast$-autonomous
structure of the category. By Lemma \ref{LemmaJoin-preserving},
the functor given in the following definition is well-defined.

\begin{definition}
Let $\F:\cCas\rightarrow\cChu_{2_0}$ be the functor defined on
objects as $\F(\l)=(\Sigma\cup\{0\},r,\Sigma'\cup\{1\})$,
with $\Sigma\cup\{0\}$ pointed by $0$ and $\Sigma'\cup\{1\}$
pointed by $1$, and with $r(p,x)=0\iff p\leq x$ for any
$p\in\Sigma\cup\{0\}$ and $x\in\Sigma'\cup\{1\}$ (hence
$r(0,\cdot)\equiv r(\cdot,1)\equiv 0$ ). On arrows, the functor
$\F$ is defined as $\F(f)=(f,f^\circ)$.
\end{definition}

\begin{lemma}\label{LemmaFFullFaithfull}
The functor $\F:\cCas\rightarrow\cChu_{2_0}$ is full and faithful.
In particular, $\F(\l_1)\cong\F(\l_2) \implies \l_1\cong\l_2$.
\end{lemma}

\begin{proof} Let $\l_1,\,\l_2\in\cCas$. To prove that $\F$ is
faithful, let $f$
and
$g\in\cCas(\l_1,\l_2)$ be such that $\F(f)=\F(g)$.
Thus, $f=g$ on atoms, and since those maps preserve arbitrary
joins, for any $a\in\l_1$ we have $f(a)=\bigvee f(\s[a])=\bigvee
g(\s[a])=g(a)$.

To show that $\F$ is full, let $(f,g):\F(\l_1)\rightarrow\F(\l_2)$
be an arrow of $\cChu_{2_0}$. Write $\F(\l_i)$ as
$(\s_i\cup\{0\},r_i,\s'_i\cup\{1\})$. Define
$h:\l_1\rightarrow\l_2$ as $h(a)=\bigvee f(\s[a])$. Since $f$ is a
pointed map from $\s_1\cup\{0_1\}$ to $\s_2\cup\{0_2\}$, we
have that $h$ preserves $0$ and sends atoms to atoms or $0$.
Denote by $H$ the restriction of $h$ to atoms (hence $H=f$). Let
$x$ be a coatom of $\l_2$. Then
\[\begin{split}H^{-1}(\s[x]\cup\{0\})&=\{p\in\s_1\,;\,r_2(f(p),x)=0\}\\
&= \{p\in\s_1\,;\,r_1(p,g(x))=0\}=\s[g(x)]. \end{split}\]
Therefore, by Lemma \ref{LemmaJoin-preserving}, $h$ preserves
arbitrary joins, and moreover, we find that
$h\in\cCas(\l_1,\l_2)$.

Finally, if $\F(\l_1)\cong \F(\l_2)$, then
$\clo(\l_1)=\clo(\l_2)$, therefore $\l_1\cong\l_2$.
\end{proof}

\begin{lemma}\label{LemmaCalSymClosedunderdual}
Let $\l,\, \l_1,\,\l_2\in\cCas$ and $f\in\cCas(\l_1,\l_2)$. We
have that $\F(\l^\op)=\F(\l)^\bot$ and $\F(f^\op)=\F(f)^\bot$, so
that $\F(\cCas)$ is closed under $-^\bot$. Moreover,
$\F(2)\cong\top$.
\end{lemma}

\begin{proof} By definition,
\[\F (\l^\op)=(\Sigma'\cup\{1\},r^\op,\Sigma\cup\{0\}),\]
with $r^\op(1,\cdot)\equiv 0\equiv r^\op(\cdot,0)$, and for any
$x\in\s'$ and $p\in\s$, $r^\op(x,p)=0\iff x\leq^\op
p\iff p\leq x$. As a consequence, $r^\op=\check{r}$ and
$\F(\l^\op)=\F(\l)^\bot$.

Moreover, $\F (f^\op)=(f^\op,(f^\op)^\circ)$. Now,
$f^\op:\l_2^\op\rightarrow\l_1^\op$, so that
$(f^\op)^\circ:\l_1^\op\rightarrow\l_2^\op$. Let $p$ be an atom of
$\l_1$. Then
\[\begin{split}(f^\op)^\circ(p)&=\bigvee\nolimits^\op\{b\in\l_2^\op\,;\,
f^\op(b)\leq^\op
p\}\\
&=\bigwedge\{b\in\l_2\,;\,p\leq
f^\circ(b)\}=\bigwedge\{b\in\l_2\,;\,f(p)\leq b\}=f(p).\end{split}\]
Recall that by definition, for any $b\in\l_2^\op$, we have
$f^\op(b)=f^\circ(b)$. As a consequence,
$\F(f^\op)=(f^\circ,f)=\F(f)^\bot$.

By definition, $\F(2)=(A,r,X)$ with $A=\{1,0\}$ pointed by $0$ and
$X=\{0,1\}$ pointed by $1$, and with $r(0,\cdot)\equiv 0\equiv
r(\cdot,1)$, and $r(1,0)=1$. As a consequence, $\F(2)\cong \top$.
\end{proof}

\begin{lemma}\label{CalSymClosedunderTensor}
Let $\cC$ be the subcategory of $\cChu_{2_0}$ formed by closing
$\F(\cCas)$ under isomorphisms.
Then $-\otchu-$ is a bifunctor in $\cC$.
Moreover, there are isomorphisms
\[\alpha_{_{\l_1\l_2}}:\F (\l_1)\otchu\F (\l_2)\rightarrow\F
(\l_1\circledast\l_2),\]
natural in $\l_1$ and $\l_2$.
\end{lemma}

\begin{proof} Let $\l_1,\,\l_2\in\cCas$. By Definition
\ref{DefinitionTensorChu}
\[\F (\l_1)\otchu\F (\l_2)=
((\Sigma_1\cup\{0\})\smash(\Sigma_2\cup\{0\}),t,\cChu_{2_0}(\F
(\l_1),\F (\l_2)^\top)),\]
and by Lemma \ref{Lemma*circcoatomistic},
\[\F (\l_1\circledast\l_2)=
((\Sigma_1\times\Sigma_2)\cup\{0\},t',\s_\circledast'\cup\{1\}),\]
where $t'(p,x)=0\iff p\in x$, for any $p\in\s_1\times\s_2$
and $x\in\s_\circledast'$.

Define
\[\chi:(\Sigma_1\cup\{0\})\smash(\Sigma_2\cup\{0\})\rightarrow
(\Sigma_1\times\Sigma_2)\cup\{0\};\, \chi(0_\sharp)=0,\,
\chi(p_1,p_2)=(p_1,p_2),\]
and let $\xi$ be the bijection of Lemma
\ref{LemmaBijectionCoatoms*TensorCalSym(L1,L2^op)} (write
$\xi(x)=f_x$). Moreover, define $\alpha:=(\chi,\F\circ\xi)$. By
definition of the smash product, the map $\chi$ is a bijection.

Write $\F(\l_i)$ as $(\s_i\cup\{0\},r_i,\s'_i\cup\{1\})$ and
let $p\in\s_1\times\s_2$ and $x\in\s_\circledast'$. Then we have
\[\begin{split}t'(\chi(p),x)=0&\iff
p_2\in x_2[p_1]\iff p_2\leq f_x(p_1)\iff
r_2(p_2,f_x(p_1))=0\\
&\iff t((p_1,p_2),\F(f_x))=0.\end{split}\]

Let $(f,g)\in\cChu_{2_0}(\F(\l_1),\F(\l_2)^\bot)$ and
$h\in\cCas(\l_1,\l_2^\op)$ with $\F(h)=(f,g)$. Then
\[\begin{split}t(\chi^{-1}(p),(f,g))=0&\iff
r_2(p_2,f(p_1))=0\\ &\iff r_2(p_2,h(p_1))=0\iff
t'(p,\xi^{-1}(h))=0.\end{split}\]
As a consequence, we find that
$\alpha:\F(\l_1)\otchu\F(\l_2)\rightarrow\F(\l_1\circledast\l_2)$
is an invertible arrow of $\cChu_{2_0}$.

Finally, we prove that
\[\F (f_1\circledast
f_2)\circ\alpha=\alpha\circ(\F (f_1)\otchu\F (f_2)).\]
Let $\l_1^2,\,\l_2^2\in\cCas$, $f_1\in\cCas(\l_1,\l_1^2)$ and
$f_2\in\cCas(\l_2,\l_2^2)$. Moreover, let
$p=(p_1,p_2)\in\s_1\times\s_2$ and $x$ be a coatom of
$\l_1^2\circledast\l_2^2$. Then
\[\bigvee ((f_1\circledast
f_2)^\circ(x))_i[p]=f^\circ_i(\bigvee x_i[f(p)])\]
(see the proof
of Lemma \ref{LemmaAertsOtfBiFuntors}). Therefore, we find that
\[\begin{split}(\F\circ \xi\circ(f_1\circledast f_2)^\circ)&(x)
=(f_2^\circ\circ f_x\circ
f_1,(f_2^\circ\circ f_x\circ f_1)^\circ)\\
&=(f_2^\circ\circ f_x\circ f_1,f_1^\circ\circ f_x^\circ\circ
f_2)=(\F(f_1)\otchu\F(f_2))((\F\circ \xi)(x)).\end{split}\]
\end{proof}

\begin{theorem}\label{TheoremTheTheorem}
The category $\langle\cCas, \circledast,2,\multimap,2\rangle$
with $\l_1\multimap\l_2={(\l_1\circledast \l_2^\op)}^\op$ is
$\ast$-autonomous.
\end{theorem}

\begin{proof} This theorem can be proved directly, however it also
follows easily from Proposition
\ref{Corollary*TensorisaBiFunctor}, the Lemmata
\ref{LemmaFFullFaithfull}, \ref{LemmaCalSymClosedunderdual},
\ref{CalSymClosedunderTensor}, and Proposition
\ref{LemmaChuPointedSetsis*autonomous}.
\end{proof}

\begin{remark}
We now give two examples where Axiom ${\mathbf A_0}$ does not hold
in $\l_1$ and $\l_2$, and where the set $\s_\circledast'$ defined
in Definition \ref{Definition*TensorProduct} is not a set of
coatoms. More precisely, $\l_1=\l_2=\l$ and there is
$R,\,S\in\s_\circledast'$ with $R$ a proper subset of $S$ and
$S\ne \s_1\times\s_2$.

Now, in the proof of Lemma \ref{CalSymClosedunderTensor}, we have
shown that
$\F(\l_1)\otchu\F(\l_2)\cong((\s_1\times\s_2)\cup\{0\},t',\s_\circledast'\cup\{1\})$.
As a consequence, since $R\subsetneqq S\subsetneqq\s_1\times\s_2$,
there is no complete atomistic coatomistic lattice $\l_0$ such
that $\F(\l_0)\cong\F(\l_1)\otchu\F(\l_2)$.

In Example \ref{ExamplePowerSets} $\l$ is a powerset lattice
whereas in Example \ref{ExampleA0Necessary} $\l$ is an irreducible
complete atomistic orthocomplemented lattice. Recall that Axiom
${\mathbf A_0}$ implies irreducibility.
\end{remark}

\begin{example}\label{ExamplePowerSets}
Let $\l=2^{\s}$ be a powerset lattice and $r_0\in\s$ (with $\#\s\geq 2$). Let
$R=\bigcup\{\{r\}\times(\s\backslash\{r\})\,;\,r\in\s\}$ and
$S=R\cup\{r_0\}\times\s$. Then $R$ and $S$ are subsets of
$\s\times\s$, and obviously, for all $p\in\s\times\s$, $R_1[p]$, $R_2[p]$,
$S_1[p]$, and $S_2[p]$ are in $\clo(\s'\cup\{1\})$. Moreover,
$R$ is a proper subset of $S$, and $S\ne\s\times\s$.
\end{example}

\begin{example}\label{ExampleA0Necessary}
Let $\l$ be the orthocomplemented simple closure space on
$\s=\mathbb{Z}_{6}$ with $n'=\{(n+2),(n+3),(n+4)\}$, where
$(m):=m\ {\rm mod}\ 6$ and $-'$ denotes the orthocomplementation.
Use the map $g:\s\rightarrow\mathbb{C};\, n\mapsto{\rm e}^{{\rm
i}n\pi/3}$ to check that $n\p m\iff n\in m'$ is indeed
symmetric, anti-reflexive and separating, {\it i.e.} for any
$p,\,q\in\s$ there is $r\in\s$ such that $p\p r$ and $q\not\p r$.
Obviously, $\l$ is irreducible, but $\l$ does not satisfy Axiom
${\mathbf A_0}$. For instance, $0'\cup 3'=\s$.

Let $R\subseteq\s\times\s$ defined as
$R=\{0,1,2\}\times\{0,1,2\}\cup\{3,4,5\}\times\{3,4,5\}$ and
$S:=4'\times\s\cup\s\times 1'$. Obviously, $R$ and $S$ are in
$\s_\circledast'$. Moreover, $R$ is a proper subset of $S$, and
$S\ne\s\times\s$.
\end{example}

\begin{remark}
Note that the category $\cCas$ is different from the $\ast$-autonomous category of
complete
semilattices
discussed in \cite{Barr:handbook}.
For instance, this category is complete whereas $\cCas$ lacks coproducts.
\end{remark}

\begin{theorem}
The category $\cCas$ does not have binary coproducts.
\end{theorem}

\begin{proof}
We proceed \textit{ad absurdum} assuming that there is a an object $\l_1+\l_2\in\cCas$ and
arrows $i_1\in\cCas(\l_1,\l_1+\l_2)$ and $i_2\in\cCas(\l_2,\l_1+\l_2)$ such that for any object $\l\in\cCas$
and any arrows $h_1\in\cCas(\l_1,\l)$ and $h_2\in\cCas(\l_2,\l)$, there is a unique arrow $h\in\cCas(\l_1+\l_2,\l)$
such that $h\circ i_1=h_1$ and $h\circ i_2=h_2$:
\begin{equation}\label{coproduct}
\bfig\Vtrianglepair/{->}`{<-}`{->}`{-->}`{->}/[\l_1`\l_1+\l_2`\l_2`\l;i_1`i_2`h_1`h`h_2]\efig
\end{equation}
We denote by $\s$ (respectively $\s'$) the set of atoms (respectively of coatoms) of $\l_1+\l_2$.
For $\l\in\cCas$ and a coatom $x$ of $\l$, we write $h_x:\l\rightarrow 2$ for the map defined by
\[
h_x(0)=0,\ h_x(p)=\left\{\begin{array}{ll}0&\mbox{if }p\in\s[x]\\[2mm]1&\mbox{if }p\in\s\backslash\s[x]\end{array}\right.
\ \mbox{and}\ h_x(a)=\bigvee h_x(\s[a])\
\forall\ a\in\l.
\]
Note that $h_x\in\cCas(\l,2)$ by Lemma \ref{LemmaJoin-preserving}. For $\mathcal{G}\in\cCas$ and the constant map $h=0$ from $\l$ to $\mathcal{G}$ one also has
$h\in\cCas(\l,\mathcal{G})$ by Lemma \ref{LemmaJoin-preserving}.

(1) {\bf Claim:} $i_1$ and $i_2$ are injective. [{\it Proof.} Put $\l=\l_1$, $h_1=\mathsf{id}$ and
$h_2=0$. Then
\[
i_1(a)=i_1(b)\implies h(i_1(a))=h(i_1(b))\implies h_1(a)=h_1(b)\implies a=b.
\]
The proof for $i_2$ is similar.]

(2)  {\bf Claim:} $i_1(\s_1)\subseteq\s$ and $i_2(\s_2)\subseteq\s$. [{\it Proof.}
Let $p_1\in\s_1$. By definition, $i_1(p_1)$ is an atom of $\l_1+\l_2$ or $0$. Suppose that $i_1(p_1)=0$. Then
\[
p_1=h_1(p_1)=h(i_1(p_1))=h(0)=0,
\]
a contradiction. The proof for $i_2$ is similar.]

(3)  {\bf Claim:} $i_1(\s_1)\cap i_2(\s_2)=\emptyset$. [{\it Proof.} Put $\l=\l_1$, $h_1=\mathsf{id}$ and
$h_2=0$. Suppose that there is $p_1\in\s_1$ and $p_2\in\s_2$ such that $i_1(p_1)=i_2(p_2)$. Then
\[
p_1=h_1(p_1)=h(i_1(p_1))=h(i_2(p_2))=h_2(p_2)=0,
\]
a contradiction.]

Define $\Omega$ by
\[
\s=i_1(\s_1)\sqcup i_2(\s_2)\sqcup\Omega.
\]

(4) {\bf Claim:} $i_1(1)\vee i_2(1)=1$.  [{\it Proof.} Suppose that there is a coatom $x\in\s'$ such that
$i_1(1)\vee i_2(1)\leq x$. Put $\l=2$ and $h_1=0=h_2$.
Let $h\in\cCas(\l_1+\l_2,2)$ be the constant map $h=0$
and let $h_x\in\cCas(\l_1+\l_2,2)$. Then both $h$ and $h_x$ make the diagram
\ref{coproduct} commute, a contradiction.]

(5) {\bf Claim:} For any $x_1\in\s_1'$ and any $x_2\in\s_2'$, $i_1(x_1)\vee i_2(1)\in\s'$ and
$i_1(1)\vee i_2(x_2)\in\s'$. [{\it Proof.} First, $i_1(x_1)\vee i_2(1)\ne 1$. Indeed, put $\l=\l_1$, $h_1=\mathsf{id}$ and
$h_2=0$. Then $h^\circ(x_1)$ is a coatom or $1$. Obviously,  $h^\circ(x_1)\ne 1$ since for any atom $q_1\in\s_1$ with
$q_1\land x_1=0$, $h(i_1(q_1))=h_1(q_1)=q_1$ is not under $x_1$. As a consequence, $h^\circ(x_1)$ is a coatom and
$h^\circ(x_1)\geq i_1(x_1)\vee i_2(1)$
since
\[
h(i_1(x_1)\vee i_1(1))=h(i_1(x_1))\vee h(i_2(1))=h_1(x_1)\vee h_2(1)=x_1\vee 0=x_1.
\]

Now, suppose that $i_1(x_1)\vee i_2(1)$ is not a coatom. Then, since $\l_1+\l_2$
is coatomistic, there is at least two coatoms, say $x$ and $y$, above $i_1(x_1)\vee i_2(1)$.
Note that
\[
\s[x]=i_1(\s_1[x_1])\sqcup i_2(\s_2)\sqcup R_x
\]
with $R_x\subseteq\Omega$. Indeed, if $x\backslash \{i_1(\s_1[x_1])\sqcup i_2(\s_2)\}\nsubseteq\Omega$,
then there is an atom $p_1\in\s_1$ with $p_1\land x_1=0$
and
$x\geq i_1(p_1)\vee i_1(x_1)\vee i_2(1)=i_1(1)\vee i_2(1)=1$, a contradiction. Similarly,
\[
\s[y]=i_1(\s_1[x_1])\sqcup i_2(\s_2)\sqcup R_y.
\]
Put $\l=2$, $h_1=h_{x_1}$ and $h_2=0$. Then $h_x\in\cCas(\l_1+\l_2,2)$
and $h_y\in\cCas(\l_1+\l_2,2)$ make the diagram
\ref{coproduct} commute, a contradiction.
As a consequence,
\[
h^\circ(x_1)=i_1(x_1)\vee i_2(1).
\]
The proof of $i_1(1)\vee i_2(x_2)\in\s'$ is similar.]

(6) {\bf Claim:} For any coatom $x\in\s'$ there is $x_1\in\s_1'$ or $x_2\in\s_2'$ such that
$x=i_1(x_1)\vee i_2(1)$ or $x=i_1(1)\vee i_2(x_2)$. [{\it Proof.} By definition,
$i_1^\circ(x)$ and $i_2^\circ(x)$ are coatoms or $1$. Moreover, $i_1(i_1^\circ(x))\vee i_2(i_2^\circ(x))\leq x$.
As a consequence, by part~4, $i_1^\circ(x)$ and $i_2^\circ(x)$ cannot both be equal to $1$. Suppose now that both,
$i_1^\circ(x)$ and $i_2^\circ(x)$ are coatoms, say $x_1$ and $x_2$. Then $i_1(x_1)\vee i_2(x_2)=x$. Indeed, otherwise,
since $\l_1+\l_2$
is coatomistic, there is at least one coatom, say $y$, different from $x$ and above $i_1(x_1)\vee i_2(x_2)$.
Again,
\[
\begin{split}
\s[x]&=i_1(\s_1[x_1])\sqcup i_2(\s_2[x_2])\sqcup R_x\\
\s[y]&=i_1(\s_1[x_1])\sqcup i_2(\s_2[x_2])\sqcup R_y.
\end{split}
\]
Put $\l=2$, $h_1=h_{x_1}$ and $h_2=h_{x_2}$. Then $h_x\in\cCas(\l_1+\l_2,2)$
and $h_y\in\cCas(\l_1+\l_2,2)$ make the diagram
\ref{coproduct} commute, a contradiction. Now, if $x=i_1(x_1)\vee i_2(x_2)$, then $x\leq i_1(x_1)\vee i_2(1)$ which is
a coatom by part 5.]

(7) {\bf Claim:} $\Omega=\emptyset$. [{\it Proof.} Suppose that $\Omega\ne \emptyset$. Let $p\in\Omega$.
Put $\l=\l_1$, $h_1=\mathsf{id}$ and
$h_2=0$. Assume that there is $p\in\Omega$ such that $h(p)\ne 0$. Write $p_1=h(p)$.
Let $x=i_1(x_1)\vee i_2(1)$ be a coatom above $p$. Then
\[
\begin{split}
p_1=h(p)\leq h(x)&=h(i_1(x_1)\vee i_2(1))\\
&=h(i_1(x_1))\vee h(i_2(1))=h_1(x_1)\vee h_2(1)=x_1\vee 0=x_1,
\end{split}
\]
so that $i_1(p_1)\leq i_1(x_1)\leq i_1(x_1)\vee i_2(1)= x$. Let $y=i_1(1)\vee i_2(x_2)$ a
coatom above $p$.
Then
\[
\begin{split}
p_1=h(p)\leq h(y)&=h(i_1(1)\vee i_2(x_2))\\
&=h(i_1(1))\vee h(i_2(x_2))=h_1(1)\vee h_2(x_2)=1\vee 0=1
\end{split}
\]
so that $i_1(p_1)\leq i_1(1)\leq i_1(1)\vee i_2(x_2)= y$.
Since $\l_1+\l_2$ is atomistic, we find that
$p=\cap\{x\in\s'\,;\,p\leq x\}$, and by what precedes,
\[
p=\cap\{x\in\s'\,;\,p\leq x\}\geq p\vee i_1(p_1),
\]
a contradiction. As a consequence, $h(\Omega)=0$. Therefore, for all coatoms $x_1\in\s_1'$,
we have by part 5
\[
\begin{split}
\s[i_1(x_1)\vee i_2(1)]=\s[h^\circ(x_1)]&=i(\s_1[x_1])\sqcup i_2(\s_2)\sqcup (\s[h^\circ(x_1)]\cap\Omega)\\
&=i(\s_1[x_1])\sqcup i_2(\s_2)\sqcup \Omega.
\end{split}
\]
Put now $\l=\l_2$, $h_1=0$ and $h_2=\mathsf{id}$. A similar argument shows that $h(\Omega)=0$.
As a consequence,
for all coatoms $x_2\in\s_2'$, we have
\[
\s[i_1(1)\vee i_2(x_2)]=i(\s_1)\sqcup i_2(\s_2[x_2])\sqcup \Omega.
\]
Therefore, since $\l_1+\l_2$ is atomistic,
\[
\emptyset=\cap\{\s[x]\,;\, x\in\s'\}=\Omega.
\]
Hence, we have proved that $\s=i_1(\s_1)\sqcup i_2(\s_2)$.]

(8) {\bf Claim:} Axiom ${\mathbf A_0}$ does not hold in $\l_1+\l_2$. [{\it Proof.} Let $x_1\in\s_1'$ and
$x_2\in\s_2'$. Then, by part 7, any atom of $\l_1+\l_2$ is under $i_1(x_1)\vee i_2(1)$ or $i_1(1)\vee i_2(x_2)$.]
\end{proof}
\section{Comparison of $\circledast$ with other tensor
products}\label{SectionComparison}

In this section, we compare the tensor product $\circledast$ with
the {\it separated product} of Aerts \cite{Aerts:1982}, the
 box product of Gr\"atzer and Wehrung \cite{Graetzer/Wehrung:1999},
and with the tensor product of Shmuely \cite{Shmuely:1974}.

\begin{definition}\label{DefinitionSeparatedProduct}
Let $\l_1$ and $\l_2$ be complete atomistic lattices. Write
$\aut(\l_i)$ for the group of automorphisms of $\l_i$. Define
\[\begin{split}
\l_1\varowedge\l_2&:=\lac\bigcap \omega\,;\, \omega\subseteq
\{a_1\mathbin{\square} a_2\,
;\, a\in\l_1\times\l_2\}\rac\, ,\\
\l_1\varovee\l_2&:=\{R\subseteq\s_1\times\s_2\,;\,
R_1[p]\in\clo(\l_1)\ \mbox{and}\ R_2[p]\in\clo(\l_2),\, \forall\,
p\in\s_1\times\s_2\},
\end{split}\]
ordered by set-inclusion.

Moreover, $\Sep(\l_1,\l_2)$ is defined as the set of all simple
closure spaces $\l$ on $\s_1\times\s_2$ such that
\begin{enumerate}
\item $\l_1\varowedge\l_2\subseteq\l\subseteq\l_1\varovee\l_2$,
and
\vspace{0.1cm}\item for all $(u_1,u_2)\in\aut(\l_1)\times\aut(\l_2)$,
there is $u\in\aut(\l)$ such that $u(p_1,p_2)\hspace{-0.1cm}=$
$(u_1(p_1),u_2(p_2))$, for all atoms $p_1\in\s_1$ and
$p_2\in\s_2$.
\end{enumerate}
\end{definition}

\begin{remark}\label{RemarkCoatomsAerts*Chu} By Lemma
\ref{LemmaCoatomsAertsCoatoms*} and from the proof of Lemma
\ref{Lemma*circcoatomistic}, it follows that if
$\l_1,\,\l_2\in\cCas$, then
$\s_\varowedge'\subseteq\s_\circledast'\subseteq\s_\varovee'$,
where $\s_\varovee'$ denotes the set of coatoms of
$\l_1\varovee\l_2$. Note that obviously, $\s_\varowedge'$ is the
set of coatoms of $\l_1\varowedge\l_2$. Moreover,
$\l_1\varowedge\l_2$ is coatomistic (see \cite{Ischi:2007}).
\end{remark}

\begin{theorem}\label{TheoremWeakTensorProducts} Let $\l_1$ and $\l_2$
be complete atomistic lattices. Then $\l_1\varowedge\l_2$ and
$\l_1\varovee\l_2$ are simple closure spaces on $\s_1\times\s_2$.
Moreover, $\Sep(\l_1,\l_2)$ (ordered by set-inclusion) is a
complete lattice.
\end{theorem}

\begin{proof} See \cite{Ischi:2007}. \end{proof}

\begin{remark}\label{RemarkSeparatedProduct}
If $\l_1$ and $\l_2$ are orthocomplemented, so is
$\l_1\varowedge\l_2$; it is the {\it separated product} of Aerts
\cite{Aerts:1982} (see \cite{Ischi:2007}). The binary relation on
$\s_1\times\s_2$ defined by $p\mathbin{\#}q\iff p_1\p_1 q_1$ or
$p_2\p_2 q_2$, induces an orthocomplementation of
$\l_1\varowedge\l_2$.

For atomistic lattices, define $\l_1\varowedge_n\l_2$ by taking
only finite intersections in Definition
\ref{DefinitionSeparatedProduct}. Then
$\l_1\varowedge_n\l_2\cong\l_1\mathbin{\square}\l_2$ which is the {\it
box-product} of Gr\"atzer and Wehrung \cite{Graetzer/Wehrung:1999}
(see \cite{Ischi:2007}).\end{remark}

By Lemma \ref{LemmaJoin-preserving}, the functor given in the
following definition is well-defined.

\begin{definition} Let $\cCa$ be the category of complete atomistic
lattices with maps preserving arbitrary joins and sending atoms to
atoms. We denote by $\G:\cCa\rightarrow\cChu(\cSet,2)$ the functor
defined on objects as $\G(\l):=(\s,r,\l)$ with
$r(p,a)=1\iff p\leq a$, and on arrows as
$\G(f)=(f,f^\circ)$.
\end{definition}

\begin{proposition}\label{PropositionFunctorG} The functor $\G$
is full and faithful. Moreover for any $\l_1,\,\l_2\in\cCa$, we
have $\G(\l_1\varovee\l_2)\cong\G(\l_1)\otchu\G(\l_2)$, where
$-\otchu-$ is the bifunctor in the category $\cChu(\cSet,2)$.
\end{proposition}

\begin{proof}
To prove that $\G$ is full and faithful, we can proceed as in the
proof of Lemma \ref{LemmaFFullFaithfull}. For the rest of the
proof, we refer to \cite{Ischi:2007}.
\end{proof}

\begin{remark}\label{RemarkShmuely} For $\l_1$ and $\l_2$
complete atomistic lattices, we have
$\l_1\varovee\l_2\cong\l_1\ot\l_2$ the tensor product of Shmuely
\cite{Shmuely:1974} (see \cite{Ischi:2007}).
\end{remark}

\begin{theorem}Let $\l_1,\,\l_2\in\cCas$. Then
$\l_1\circledast\l_2\in\Sep(\l_1,\l_2)$. \end{theorem}

\begin{proof} By definition \ref{Definition*TensorProduct},
$\s_\circledast'\subseteq\l_1\varovee\l_2$. Now, by Theorem
\ref{TheoremWeakTensorProducts}, $\l_1\varovee\l_2$ is a simple
closure space on $\s_1\times\s_2$. As a consequence, from Definition
\ref{Definition*TensorProduct} it follows that
$\l_1\circledast\l_2\subseteq\l_1\varovee\l_2$.

Similarly, by Lemma \ref{LemmaCoatomsAertsCoatoms*},
$\s_\varowedge'\subseteq\s_\circledast'$, hence by Remark
\ref{RemarkCoatomsAerts*Chu} and Definitions
\ref{Definition*TensorProduct} and
\ref{DefinitionSeparatedProduct}, it follows that
$\l_1\varowedge\l_2\subseteq\l_1\circledast\l_2$.

Finally, by Lemma \ref{LemmaAertsOtfBiFuntors}, for all
automorphisms $u_1$ and $u_2$, $u_1\circledast
u_2\in\aut(\l_1\circledast\l_2)$.
\end{proof}

We end this section with two examples.

\begin{example}
As in Example \ref{ExampleA0Necessary}, let $\l$ be the
orthocomplemented simple closure space on $\s=\mathbb{Z}_{12}$
with $n'=\{(n+5),(n+6),(n+7)\}$, where $(m):=m\ {\rm mod}\ 12$.
Note that Axiom ${\mathbf A_0}$ holds in $\l$ but that $\l$ does
not have the covering property since $2\vee 8=\s\gneqq
2'\gvertneqq 8$.

Let
\[x:=(0'\times 0')\cup (3'\times 3')\cup (6'\times 6')\cup (9'\times 9').\]
By Definition \ref{Definition*TensorProduct},
$x\in\l\circledast\l$ (note that $x\in\s_\circledast'$), whereas
obviously, $x^\#=\{p\in\s\times\s\,;\,p\#q\,\forall q\in x\}=\emptyset$ (see Remark
\ref{RemarkSeparatedProduct}), therefore $x\not\in\l\varowedge\l$.

Let
\[R:= (0'\times 0')\cup((2'\cap 3')\times (2'\cap 3'))\cup (5'\times 5')\cup
(8'\times 8')\cup (4\times 4).\]
By Definition \ref{DefinitionSeparatedProduct},
$R\in\l\varovee\l$. {\bf Claim:} $R$ is a coatom of
$\l\varovee\l$. [{\it Proof.}  Let $p=(p_1,p_2)$ be an atom not
under $R$. Define $z:=p\vee R$, where the join is taken in
$\l\varovee\l$. If $p_1$ or $p_2$ is in $0'$, $5'$ or $8'$, then
$p\vee R=\s\times\s$. Indeed, suppose for instance that $p_1$ is in
$0'$. Then $p_1\times ( 0'\cup\{p_1\})\subseteq z$, hence
$p_1\times\s\subseteq z$. Therefore,
\[(p_1\vee 5')\times 5'\cup (p_1\vee 8')\times
8'\subseteq z.\]
Now, since $0'$, $5'$ and $8'$ are disjoint, $p_1\not\in 5'$ and
$p_1\not\in 8'$, thus $p_1\vee 5'=\s$ and $p_1\vee8'=\s$. As
a consequence, $\s\times( 5'\cup 8')\subseteq z$, whence
$z=\s\times\s$.

Finally, suppose for instance that $p_1=4$ and $p_2\in 2'\cap 3'=\{8,9\}$. Then $4\times(4\vee p_2)=4\times \s\subseteq
z$, hence
\[(4\vee 0')\times 0'\cup (4\vee 5')\times
5'\cup (4\vee 8')\times 8'\subseteq z,\]
therefore $\s\times(0'\cup 5'\cup 8')\subseteq z$. As a
consequence, $z=\s\times\s$.]

Obviously, $R_2[(4,\cdot)]=\{4\}\not\in\s'$, hence $R$ is not a
coatom of $\l\circledast\l$. As a consequence,
$\l\circledast\l\ne\l\varovee\l$ (see Remark
\ref{RemarkCoatomsAerts*Chu}).

To summarize, we have $\l\varowedge\l\varsubsetneqq
\l\circledast\l\varsubsetneqq\l\varovee\l$.
\end{example}

\begin{example}\label{Example3tensor=aertsouChu}
We leave it as an exercise to prove that
\[
\mathsf{MO}_3\circledast\mathsf{MO}_4=
\mathsf{MO}_3\,\varowedge\mathsf{MO}_4,\]
where $\mathsf{MO}_n$ was defined in Section \ref{setting}.
\end{example}
\section{Weak bimorphisms}\label{SectionWeakBimorphisms}

We now prove that the tensor product $\circledast$ can be defined
as the solution of a universal problem with respect to {\it weak
bimorphisms}.

\begin{definition}\label{DefinitionWeakBimorphisms}
Let $\l_1,\,\l_2,\,\l\in\cCas$ and $f:\l_1\times\l_2\rightarrow\l$ be
a map. Then $f$ is a {\it weak bimorphism} if for any $p_1\in\s_1$
and for any $p_2\in\s_2$, we have $f(-,p_2)\in\cCas(\l_1,\l)$ and
$f(p_1,-)\in\cCas(\l_2,\l)$. Moreover, $\l$ is a {\it $w-$tensor
product} of $\l_1$ and $\l_2$ if there is a weak bimorphism
$f:\l_1\times\l_2\rightarrow\l$ such that for any $\l_0\in \cCas$ and
any weak bimorphism $g:\l_1\times\l_2\rightarrow \l_0$, there is a
unique arrow $h\in\cCas(\l,\l_0)$ such that the following diagram
commutes:
$$\bfig
 \ptriangle|alm|/>`>`>/[\l_1\times\l_2`\l`\l_0;f`g`!h]
 \efig$$
\end{definition}

\begin{remark}
By definition, the $w-$tensor product is unique up to
isomorphisms.
\end{remark}

\begin{theorem}\label{TheoremstarcircWeakTensorProduct}
Let $\l_1,\,\l_2\in\cCas$. Then $\l_1\circledast\l_2$ is the
$w-$tensor product of $\l_1$ and $\l_2$.
\end{theorem}

\begin{proof}
Define $f:\l_1\times\l_2\rightarrow\l_1\circledast\l_2$ as
$f(a)=\bigvee (a_1\circ a_2)$ where the join is taken in
$\l_1\circledast\l_2$. Note that
\[ a_1\circ a_2=(a_1\circ 1)\cap(1\circ a_2)=(a_1 \mathbin{\square}
0)\cap(0\mathbin{\square} a_2).\]
Hence, by Lemma \ref{LemmaCoatomsAertsCoatoms*}, $f(a)=a_1\circ
a_2$.

Let $p_1\in\s_1$ and $\omega\subseteq\l_2$. Then, obviously,
$f(p_1,x)\subseteq f(p_1,\bigvee\omega)$ for all $x\in\omega$,
hence $\bigvee\{ f(p_1,x)\,;\, x\in\omega\}\subseteq
f(p_1,\bigvee\omega)$. As a consequence, there is $B\subseteq\s_2$
such that $\bigvee\{ f(p_1,x)\,;\, x\in\omega\}=p_1\times B$ with
$B\subseteq \s[\bigvee\omega]$ and $\s[x]\subseteq B$ for all
$x\in\omega$. Now, since $\l_1\circledast\l_2\subseteq
\l_1\varovee\l_2$, $B\in\clo(\l_2)$. Therefore,
$B=\s[\bigvee\omega]$. As a consequence, $f$ is a weak bimorphism.

Let $\l_0\in\cCas$ and let $g:\l_1\times\l_2\rightarrow\l_0$ be a weak
bimorphism. Define $h:\l_1\circledast\l_2\rightarrow\l_0$ as
$h(a):=\bigvee\{g(p)\,;\, p\in a\}$. For $p\in\s_1\times\s_2$, define
$g_{p_1}:=g(p_1,-)$ and $g_{p_2}:=g(-,p_2)$. Let $x\in\s_0'$ the
set of coatoms of $\l_0$. Denote by $H$, $G_{p_1}$ and $G_{p_2}$
the restrictions to atoms of $h$, $g_{p_1}$ and $g_{p_2}$
respectively. Then
\[H^{-1}(\s[x]\cup\{0\})=\bigcup_{p_1\in\s_1} p_1\times
G_{p_1}^{-1}(\s[x]\cup\{0\})= \bigcup_{p_2\in\s_2}
G_{p_2}^{-1}(\s[x]\cup\{0\})\times p_2.\]
Now, since $g$ is a weak bimorphism,
$\bigvee G_{p_i}^{-1}(\s[x]\cup\{0\})$ is a coatom or $1$, therefore
$H^{-1}(\s[x]\cup\{0\})\in\s_\circledast'\cup\{1\}$. As a
consequence, by Lemma \ref{LemmaJoin-preserving} and Lemma
\ref{Lemma*circcoatomistic},
$h\in\cCas(\l_1\circledast\l_2,\l_0)$. Let
$h'\in\cCas(\l_1\circledast\l_2,\l_0)$ such that $h'\circ f=g$.
Then $h'$ equals $h$ on atoms, therefore $h'=h$.
\end{proof}

Note that Lemma \ref{LemmaAertsOtfBiFuntors} may be proved
directly by using this theorem.

\begin{remark}\label{RemarkFraser}
In a category $\cC$ concrete over $\cSet$, define bimorphisms as
maps $f:\A\times\B\rightarrow\C$ such that $f(a,-)\in\cC(\B,\C)$ and
$f(-,b)\in\cC(\A,\C)$ for all $a\in\A$ and $b\in\B$. Moreover,
define a tensor product as in Definition
\ref{DefinitionWeakBimorphisms} with bimorphisms instead of weak
bimorphisms. Then, for the category of join semilattices with maps
preserving all finite joins, the definition of a tensor product is
equivalent to the definition of the semilattice tensor product
given by Fraser in \cite{Fraser:1976} (note that $f(\l_1\times\l_2)$
generates $\l$ if and only if the arrow $h$ is unique). Note also
that for the category of complete atomistic lattices with maps
preserving arbitrary joins, the tensor product is given by
$\varovee$ (the proof is similar to the proof of Theorem
\ref{TheoremstarcircWeakTensorProduct}, see \cite{Ischi:2007}).
\end{remark}
\section{The $\circledast$ and $\multimap$ products for
DAC-lattices}\label{SectionDAC-lattices}

In quantum theory, the propositions describing a quantum entity are modeled by
the lattice of closed subspaces of a complex Hilbert space---which is a DAC-lattice.
The one-dimensional subspaces are the atoms and correspond to the possible physical states
of the system. The lattice of propositions of a compound system is the lattice of closed subspaces
of the tensor product of Hilbert spaces and possible states are the one-dimensional subspaces of the tensor product.
They are either product states
$\phi_1\otimes\phi_2$ or entangled states $\phi_1\otimes\phi_2+\psi_1\otimes\psi_2+\cdots$.
In the case where two systems are represented by DAC-lattices $\l_1$ and $\l_2$,
their tensor product $\l_1\circledast\l_2$ is a candidate for the description of
``separated'' systems, that is, a compound system in which the state space contains only
product states (see \cite{Aerts:1982}). We notice \textit{en passant} that
the lattice $\l_1\multimap\l_2$ appearing in the Chu construction contains all
entangled states (Remark \ref{rem:LastRemark}).

In this section, we study our new tensor product $\l_1\circledast\l_2$ in the case where $\l_1$ and $\l_2$
are DAC-lattices. More precisely,
we characterize
the coatoms of $\l_1\circledast\l_2$ in terms of semilinear maps (Theorem \ref{TheoremCoatoms*ForDacLattices}) and
compare the set of atoms of $\l_1\multimap\l_2$ with the set of one dimensional
subspaces of $E_1\ot E_2$, where $E_i$ denote the respective underlying vectors spaces (Theorem
\ref{TheoremLast}).

\begin{lemma}Let $\l_1,\,\l_2\in\cCas$ be DAC-lattices, $x_1,\,y_1\in\s_1'$, and
$x_2,\,y_2\in\s_2'$, and let $h:\s'[x_1\wedge
y_1]\rightarrow\s'[x_2\wedge y_2]$ be a bijection. Then
$x^h:=\bigcup\{z\circ h(z)\,;\,z\in\s'[x_1\wedge y_1]\}$ is a
coatom of $\l_1\circledast\l_2$, which we call a {\it
$\ast$-coatom}.
\end{lemma}

\begin{proof}
Let $p\in\s_1\times\s_2$. Since $\l_i$ and $\l_i^\op$ have the
covering property, either $p_i\leq x_i\wedge y_i$ or there is a
unique $z_i\in\s'[x_i\wedge y_i]$ such that $p_i\leq z_i$.
Therefore, either $x^h_2[p]=\s_2$ or $x^h_2[p]=\s[h(z_1)]$ for
some $z_1\in\s'[x_1\wedge y_1]$, and either $x^h_1[p]=\s_1$ or
$x^h_1[p]=\s[h^{-1}(z_2)]$ for some $z_2\in\s'[x_2\wedge y_2]$.
As a consequence $x^h\in\s'_\circledast$.
\end{proof}

\begin{lemma}\label{Lemma*Coatoms}
Let $\l_1,\,\l_2\in\cCas$ be DAC-lattices and
$f\in\cCas(\l_1,\l_2^\op)$ with $f(\l_1)$ of length $2$. Let $\xi$
be the bijection of Lemma
\ref{LemmaBijectionCoatoms*TensorCalSym(L1,L2^op)}. Then
$\xi^{-1}(f)$ is a
$\ast$-coatom.
\end{lemma}

\begin{proof}
Since by hypothesis $f(\l_1)$ is of length $2$, there exists
$x_2,\, y_2\in\s_2'$ such that
\[f(\l_1)=[0^\op,x_2\textstyle{\vee}^\op y_2],\]
where $[0^\op,x_2\vee^\op y_2]$ denotes the interval
$\{a\in\l_2^\op\,;\,0^\op\leq^\op a\leq^\op x_2\vee^\op y_2\}$.
We write $X:=\xi^{-1}(f)$. Hence, by hypothesis, we have
\[1\circ(x_2\textstyle{\wedge} y_2)\subseteq X:=\xi^{-1}(f).\]

(1) {\bf Claim:}  $\forall z\in\s_2'$, $1\circ z\nsubseteq
X$. [{\it Proof.}  If $1\circ z\subseteq X$ for some $z\in\s_2'$,
then $f(\l_1)=[0^\op,z]$, a contradiction, since by hypothesis,
$f(\l_1)$ is of length $2$.]

(2) {\bf Claim:}  $\forall x\in\s_1'$, $x\circ 1\nsubseteq
X$. [{\it Proof.}  Suppose that $x\circ 1\subseteq X$ for some
$x\in\s_1'$. Let $p$ be an atom of $\l_1$ not under $x$. Then $X$
contains $(\s[x]\cup p)\times \s[f(p)]$, hence, since
$\l_1\circledast\l_2\subseteq\l_1\varovee\l_2$, it follows that
$1\circ f(p)\subseteq X$; whence a contradiction by part 1.]

(3) {\bf Claim:}  $\forall q\in\s_2$,
$f^\circ(q)\circ((x_2\wedge y_2)\vee q)\subseteq X$. [{\it
Proof}. Let $q\in\s_2$. For any atom $p$ under $f^\circ(q)$, we
have $q\leq f(p)$, hence, since $1\circ(x_2\wedge y_2)\subseteq
X$, $p\times(\s[x_2\wedge y_2]\cup q)\subseteq X$. As a
consequence, $p\circ((x_2\wedge y_2)\vee q)\subseteq X$,
since $\l_1\circledast\l_2\subseteq \l_1\varovee\l_2$.]

(4) {\bf Claim:}  Let $q\in\s_2$ with $q\wedge
x_2\wedge y_2=0$. Then $f^\circ(q)\ne 1$. [{\it Proof.}
Suppose that $f^\circ(q)=1$. Then, by part 3,
$1\circ((x_2\wedge y_2)\vee q)\subseteq X$. Thus,
$(x_2\wedge y_2)\vee q\ne 1$. Moreover, since $\l_2^\op$ has
the covering property, $(x_2\wedge y_2)\vee q$ is a coatom;
whence a contradiction by part 1.]

(5) Let $z$ be a coatom of $\l_2$ above $x_2\wedge y_2$. {\bf
Claim:} For any atoms $p$ and $q$ of $\l_2$ with $p,\,
q\leq z$ and $p\wedge x_2\wedge y_2=0=q\wedge
x_2\wedge y_2$, we have $f^\circ(p)=f^\circ(q)$. [{\it Proof.}
Suppose that $f^\circ(p)\ne f^\circ(q)$. Now, $(x_2\wedge
y_2)\vee p=z=(x_2\wedge y_2)\vee q$. Whence, by part 3,
$(f^\circ(p)\vee f^\circ(q))\circ z=1\circ z\subseteq X$, a
contradiction by part 1.]

As a consequence, we can define a map $k:\s'[x_2\wedge
y_2]\rightarrow \s_1'$ as $k(z):=f^\circ(q)$ for any atom $q\leq
z$ such that $q\wedge x_2\wedge y_2=0$. Note that
$X=\{f^\circ(q)\circ q\,;\,q\in\s_2\}$. Indeed, by Lemma
\ref{LemmaBijectionCoatoms*TensorCalSym(L1,L2^op)}, $X=\{p\circ
f(p)\,;\,p\in\s_1\}$. Now, $(r,s)\in X\iff
s\in\s[f(r)]\iff r\in\s[f^\circ(s)]$. Hence, by what
precedes, we have
\[X=1\circ(x_2\textstyle{\wedge} y_2)\cup\{k(z)\circ
z\,;\,z\in\s'[x_2\textstyle{\wedge}
y_2]\}.\]

(6) {\bf Claim:}  $\mathrm{Im}(k)\subseteq
\s'[k(x_2)\wedge k(y_2)]$. [{\it Proof.}  First, note that
\[\s[k(x_2)\textstyle{\wedge} k(y_2)]\times(\s[x_2]\cup\s[y_2])
\subseteq X,\]
hence $k(x_2)\wedge k(y_2)\circ 1\subseteq X$. Suppose that
there is $z\in\s'[x_2\wedge y_2]$ such that $k(z)\ngeq
k(x_2)\wedge k(y_2)$. Then $k(z)\wedge k(x_2)\ne
k(z)\wedge k(y_2)$, and $(k(z)\wedge k(x_2))\circ 1\subseteq
X$ and $k(z)\wedge k(y_2)\circ 1\subseteq X$. Now, since
$\l_1^\op$ has the covering property,
\[(k(z)\textstyle{\wedge} k(x_2))\textstyle{\vee}(k(z)
\textstyle{\wedge} k(y_2))=k(z),\]
therefore $k(z)\circ 1\subseteq X$, a contradiction.]

(7) Let $p\in\s_1$ with $p\wedge k(x_2)\wedge k(y_2)=0$.
{\bf Claim:}  $f(p)\ne 1$. [{\it Proof.}  If $f(p)=1$, then
$(p\cup\s[k(x_2)\wedge k(y_2)])\times\s_2\subseteq X$, hence
$(p\vee (k(x_2)\wedge k(y_2)))\circ 1\subseteq X$. Now,
$p\vee(k(x_2)\wedge k(y_2))$ is a coatom, whence a
contradiction by part 2.]

(8) {\bf Claim:}  The map $k$ is surjective. [{\it Proof.}
Let $z$ be a coatom of $\l_1$ above $k(x_2)\wedge k(y_2)$, and
let $p$ be an atom of $\l_1$ under $z$ such that $p\wedge
k(x_2)\wedge k(y_2)=0$. By part 7, $f(p)$ is a coatom, and
since $p\vee(k(x_2)\wedge k(y_2))=z$, we find that
$k(f(p))=z$.]

(9) {\bf Claim:}  The map $k$ is injective. [{\it Proof.}
Let $t$ and $z$ be two coatoms above $x_2\wedge y_2$. Suppose
that $k(z)=k(t)$. Then $\s[k(z)]\times(\s[z]\cup\s[t])\subseteq X$,
hence $k(z)\circ 1\subseteq X$, a contradiction by part 2.]
\end{proof}

\begin{theorem}[Faure and Fr\"olicher,
\cite{Faure/Froelicher:handbook} Theorem 10.1.3] For $i=1$ and
$i=2$, let $E_i$ be a vector space over a division ring
$\mathbb{K}_i$, and $\pro(E_i)$ the lattice of all subspaces of
$E_i$. If $g:\pro(E_1)\rightarrow\pro(E_2)$ preserves arbitrary
joins, sends atoms to atoms or $0$, and if $g(\pro(E_1))$ is of
length $\geq 3$, then $g$ is induced by a semilinear map
$f:E_1\rightarrow E_2$ ({\it i.e.}
$g(\mathbb{K}v)=\mathbb{K}f(v)$, $\forall v\in E_1$).
\end{theorem}

\begin{corollary}\label{CorollaryFaure} For $i=1$ and
$i=2$, let $(E_i,F_i)$ be pairs of dual spaces, and let
$g:\l_{F_1}(E_1)\rightarrow\l_{F_2}(E_2)$ be a join-preserving
map, sending atoms to atoms or $0$ with $g(\l_{F_1}(E_1))$ of
length $\geq 3$. Then there is a semilinear map $f:E_1\rightarrow
E_2$ that induces $g$.\end{corollary}

\begin{proof} Define $h:\pro(E_1)\rightarrow\pro(E_2)$ as $h(V)=\bigvee
g(\s[V])$ where the join is taken in $\pro(E_2)$. Note that on
atoms $h=g$.

Denote by $H$ the restriction of $h$ to atoms and $G$ the
restriction of $g$ to atoms (hence $G=H$). Let $W$ be a subspace
of $E_2$ and let $p,\,q\in H^{-1}(W)$. Then $h(p)\vee h(q)$
(where the join is taken in $\pro(E_2)$) is a 2-dimensional
subspace of $E_2$, hence $h(p)\vee h(q)\in \l_{F_2}(E_2)$ (see
Remark \ref{RemarkFiniteSubspaces}). Therefore, by Lemma
\ref{LemmaJoin-preserving}, $G^{-1}(h(p)\vee
h(q))\in\clo(\l_{F_1}(E_1))$, hence $H^{-1}(h(p)\vee
h(q))\in \clo(\pro(E_1))$, therefore $\s[p\vee
q]\subseteq H^{-1}(h(p)\vee h(q))$. Moreover, from $h(p)\vee
h(q)\subseteq W$ it follows that $H^{-1}(h(p)\vee
h(q))\subseteq H^{-1}(W)$. As a consequence,
\[\s[p\textstyle{\vee} q]\subseteq H^{-1}(h(p)\textstyle{\vee} h(q))\subseteq
H^{-1}(W),\]
hence we have proved that $H^{-1}(W)\in\pro(E_1)$. Therefore, it
follows from Lemma \ref{LemmaJoin-preserving} that $h$ preserves
arbitrary joins.

As a consequence, there exists a semilinear map $f:E_1\rightarrow
E_2$ that induces $h$, hence also $g$ since $h$ equals $g$ on
atoms.\end{proof}

\begin{theorem}\label{TheoremCoatoms*ForDacLattices}
Let $\l_1,\,\l_2\in\cCas$ be DAC-lattices of length $\geq 4$, and
$X$ a coatom of $\l_1\circledast\l_2$. Let $(E_i,F_i)$ be pairs of
dual spaces such that $\l_i\cong \l_{F_i}(E_i)$ (see Theorem
\ref{TheoremMaeda}). Let $\xi$ be the bijection of Lemma
\ref{LemmaBijectionCoatoms*TensorCalSym(L1,L2^op)}. Then
$X$ is a
$\ast$-coatom, or there is a semilinear map $g$ from
$E_1$ to $F_2$ that induces $\xi(X)$. Moreover, $\xi(X)(\l_1)$ is of
length $1$, then $X$ is a coatom of $\l_1\varowedge\l_2$.
\end{theorem}

\begin{proof} First note that
$\l_2^\op\cong\l_{E_2}(F_2)$. Write $f=\xi(X)$. From Corollary
\ref{CorollaryFaure}, if $f(\l_1)$ is of length $\geq 3$, there is
a semilinear map $g$ from $E_1$ to $F_2$ that induces $f$, whereas
by Lemma \ref{Lemma*Coatoms}, if $f(\l_1)$ is of length $2$, $X$
is a
$\ast$-coatom.

Finally, if $f(\l_1)=[0^\op,x_2]$ for some coatom $x_2$ of $\l_2$,
then $1\circ x_2\subseteq X$. Now, let $q\in\s_2$ with $q\wedge
x_2=0$. Then $f^\circ(q)$ is a coatom of $\l_1$. Moreover,
$\s[f^\circ(q)]\times(x_2\cup q)\subseteq X$, hence
$f^\circ(q)\circ 1\subseteq X$. As a consequence,
$f^\circ(q)\mathbin{\square} x_2\subseteq X$, hence by Lemma
\ref{LemmaCoatomsAertsCoatoms*}, $X=f^\circ(q)\mathbin{\square} x_2$,
therefore $X\in\s_\varowedge'$.
\end{proof}

\begin{theorem}\label{TheoremLast} Let $E_1$ and $E_2$ be vector spaces of
dimension $n$. Then there is an injective map from the set
$\mathbb{P}(E_1\ot E_2)$ of one-dimensional subspaces of $E_1\ot
E_2$ to the set of atoms of $\pro(E_1)\multimap
\pro(E_2)$.\end{theorem}

\begin{proof}
Write $\l_1:=\pro(E_1)$ and $\l_2:=\pro(E_2)$. First, there is a
bijection from $E_1\ot E_2$ to the set of linear maps between
$E_1$ and $E_2$. Since both $E_1$ and $E_2$ are of dimension $n$,
any linear map induces an arrow in $\cCas(\l_1,\l_2)$. As a
consequence, there is an injective map from $\mathbb{P}(E_1\ot
E_2)$ to $\cCas(\l_1,\l_2)$.

Now, by definition $\l_1\multimap
\l_2={(\l_1\circledast\l_2^\op)}^\op$, and by Lemma
\ref{LemmaBijectionCoatoms*TensorCalSym(L1,L2^op)}, there is a
bijection between the set of coatoms of $\l_1\circledast\l_2^\op$
and $\cCas(\l_1,\l_2)\backslash\{f_0\}$, where $f_0$ denotes the
constant arrow which sends every atom of $\l_1$ to $0$. As a
consequence, there is a bijection between
$\cCas(\l_1,\l_2)\backslash\{f_0\}$ and the set of atoms of
$\l_1\multimap \l_2$.
\end{proof}

\begin{remark}\label{rem:LastRemark}
As a consequence of Theorem \ref{TheoremLast}, we find that
the set of atoms of
$\pro(E_1)\multimap
\pro(E_2)$ contains all the states of a compound system.
\end{remark}
\section{Acknowledgments}
The first author wishes to thank R. Seely, I. Ivanov and M. Barr
for useful discussions that led to this article.
\appendix
\section{The coherence conditions}

The category $\langle\cC,\ot,\top,\alpha,l,r\rangle$ is monoidal
if $l_{_{\top}}=r_{_{\top}}$, and for any objects $\A,\,\B,\,\C$
and $\D$, the diagrams
$$\bfig\Vtriangle[(\top\ot\A)\ot\B`\top\ot(\A\ot\B)`\A\ot\B;
\alpha_{_{\top\A\B}}`l_{_{\A}}\ot\mathsf{id}`l_{_{\A\ot\B}}]\efig
\hspace{1cm}
\bfig\Vtriangle[(\A\ot\top)\ot\B`\A\ot(\top\ot\B)`\A\ot\B;
\alpha_{_{\A\top\B}}`r_{_{\A}}\ot\mathsf{id}`\mathsf{id}\ot
l_{_{\B}}]\efig$$
$$\bfig\Vtriangle[(\A\ot\B)\ot\top`\A\ot(\B\ot\top)`\A\ot\B;
\alpha_{_{\A\B\top}}`r_{_{\A\ot\B}}`\mathsf{id}\ot
r_{_{\B}}]\efig$$
and
$$\bfig\BarrSquare/`>`<-`>/<2500,500>[((\A\ot\B)\ot\C)\ot\D`\A\ot(\B\ot(\C\ot\D))`
(\A\ot(\B\ot\C))\ot\D`\A\ot((\B\ot\C)\ot\D);
`\alpha_{_{\A\B\C}}\ot\mathsf{id}`\mathsf{id}\ot\alpha_{_{\B\C\D}}`\alpha_{_{\A\B\ot\C\D}}]
\morphism(400,500)<850,0>[`(\A\ot\B)\ot(\C\ot\D);\alpha_{_{\A\ot\B\C\D}}]
\morphism(1650,500)<400,0>[`;\alpha_{_{\A\B\C\ot\D}}] \efig$$
commute (see \cite{Eilenberg/Kelly:handbook}, p. 472, or
\cite{McLane:handbook}, \S VII.1). Further,
$\langle\cC,\ot,\top,\alpha,l,r,s\rangle$ is symmetric if the
diagrams
$$\bfig
\qtriangle[\A\ot\B`\B\ot\A`\A\ot\B;s_{_{\A\B}}`\mathsf{id}`s_{_{\B\A}}]
\efig
\hspace{1cm}\bfig
\Vtriangle[\B\ot\top`\top\ot\B`B;s_{_{\B\top}}`r_{_{\B}}`l_{_{\B}}]\efig
$$
and
$$\bfig\BarrSquare/`>`>`/<2300,500>[(\A\ot\B)\ot\C`(\B\ot\C)\ot\A`
(\B\ot\A)\ot\C`\B\ot(\C\ot\A);`s_{_{\A\B}}\ot\mathsf{id}`\alpha_{_{\B\C\A}}`]
\morphism(350,500)<850,0>[`\A\ot(\B\ot\C);\alpha_{_{\A\B\C}}]
\morphism(1520,500)<480,0>[`;s_{_{\A\B\ot\C}}]
\morphism(350,0)<850,0>[`\B\ot(\A\ot\C);\alpha_{_{\B\A\C}}]
\morphism(1520,0)<480,0>[`;\mathsf{id}\ot s_{_{\A\C}}] \efig$$
commute (see \cite{McLane:handbook}, \S VII.7).
\bibliographystyle{abbrv}


\end{document}